**Michiel Hazewinkel**                              1                              Burg. s'Jacob laan  18
Tel: +31-35-6937033                                                                1401BR   Bussum
Fax: +31-35-6917401                                                               The Netherlands
E-mail: michiel.hazewinkel@xs4all.nl

# Niceness theorems [1]

by

*Michiel Hazewinkel*
*Burg. s'Jacob laan 18*
*1401BR   Bussum*
*The Netherlands*
<michiel.hazewinkel@xs4all.nl>

## 1. Introduction and statement of the problems.

In this lecture I aim to raise a new kind of question. It appears that many important mathematical objects (including counterexamples) are unreasonably nice, beautiful and elegant. They tend to have (many) more (nice) properties and extra bits of structure than one would a priori expect.

The question is why this happens and whether this can be understood [2].

These ruminations started with the observation that it is difficult for, say, an arbitrary algebra to carry additional compatible structure. To do so it must be nice, i.e., as an algebra be regular (not in the technical sense of this word), homogeneous, everywhere the same, ... . It is for instance very difficult to construct an object that has addition, multiplication and exponentiation, all compatible in the expected ways.

The present scribblings are just a first attempt to identify and describe the phenomenon. Basically this is a prepreprint and it touches just the fringes of the subject. There is much more to be said and there are many more examples than remarked upon here.

---

[1] The first time I lectured on this at the WCAA conference, Wilfrid Laurier University in Waterloo, Canada, May 2008, the chairman summed up my lecture as follows: "If it is true it is beautiful, if it is beautiful it is probably true". I also lectured on the same subject at the Abel symposium meeting in Tromsø, Norway in June 2008.

The present screed expands on those first lectures a great deal. Yet, in spite of its length it is just a beginning: a first scratching at the edges of a great and fascinating problem that deserves devoted attention.

[2] There is of course the "anthropomorphic principle" answer, much like the question of the existence of (intelligent) life in this universe. It goes something like this. If these objects weren't nice and regular we would not be able to understand and describe them; we can see/understand only the elegant and beautiful ones. I do not consider this answer good enough though there is something in it. So the search is also on for ugly brutes of mathematical objects.

Also this anthropomorphic argument raises the subsidiary question of why we can only understand/describe beautiful/regular things. There are aspects of (Kolmogorov) complexity and information theory involved here.

This lecture is about lots of examples of this phenomenon such as Daniel Kan's observation that a group carries a comonoid structure in the category of groups if and only if it is a free group, the Milnor-Moore and Leray theorems in the theory of Hopf algebras, Grassmann manifolds and classifying spaces, and especially the star example: the ring of commutative polynomials over the integers in countably infinite indeterminates. This last one occurs all over the place in mathematics and has more compatible structures that can be believed. For instance it occurs as the algebra of symmetric functions in infinitely many variables, as the cohomology and homology of the classifying space **BU**, as the sum of the representation rings of the symmetric groups, as the free lambda-ring on one variable, as the representing ring of the Witt vectors, as the ring of rational representation of $GL$, as the underlying ring of the universal formal group, ... .

To start with, here is a preliminary list of the kind of phenomena I have in mind.

- A. Objects with a great deal of compatible structure tend to have a nice regular underlying structure and/or additional nice properties: "Extra structure simplifies the underlying object". As indicated above this sort of thing was the starting point.

- B. *Universal objects.* That is mathematical objects which satisfy a universality property. They tend to have:

a) a nice regular underlying structure

b) additional universal properties (sometimes seemingly completely unrelated to the defining universal property)

- C. Nice objects tend to be large and inversely large objects of one kind or another tend to have additional nice properties. For instance, large projective modules are free (Hyman Bass, [23 Bass]).

- D. Extremal objects tend to be nice and regular. (The symmetry of a problem tends to survive in its extremal solutions is one of the aspects of this phenomenon; even when (if properly looked at) there is bifurcation (symmetry breaking) going on.)

- E. Uniqueness theorems and rigidity theorems often yield nice objects (and inversely). They tend to be unreasonably well behaved. I.e. if one asks for an object with such and such properties and the answer is unique the object involved tends to be very regular. This is not unrelated to D.

Concrete examples of all these kinds of phenomena will be given below (section 2) as well as a (pitiful) few first explanatory general theorems (section 3).

The "niceness phenomenon" is not limited to theorems saying that e.g. in suitable circumstances an object is free; it also extends to counter examples: many of them are very regular in their construction. This can, for instance, take the form of a simple construction repeated indefinitely. Some examples are in section 2F below.

All in all I detect in present day mathematics a strong tendency towards the study of things that in some sense have low Kolmogorov complexity.

## 2. Examples

### 2A. Lots of compatible structure examples.

2A.1. *Groups in the category of groups*. To start with  here is an observation of Daniel Kan, [99 Kan] which has moreover the distinction  of being one of the first results of this kind and of admitting a nice (sic!) pictorial illustration.

First, here is the abstract setting. Let  $\mathcal{C}$   be a category with a terminal object and products. For example the category **Group**  of groups where the product is the direct product and the terminal object is the one element group.

A group object in such a category $\mathcal{C}$  is an object  $G \in \mathcal{C}$  equipped with a morphism  $m: G \times G \longrightarrow G$  (multiplication), a morphism  $e : T \longrightarrow G$  (unit element) where  $T$  is the terminal object of the category  $\mathcal{C}$ , and a morphism  $\iota : G \longrightarrow G$  (inverse)  such that the categorical versions of the standard group axioms hold. This means that the following diagrams are supposed to be commutative.

$$
\begin{array}{ccc}
G \times G \times G & \xrightarrow{\ m \times \mathrm{id}\ } & G \times G \\
{\scriptstyle \mathrm{id}\, \times\, m}\big\downarrow & & \big\downarrow{\scriptstyle m} \\
G \times G & \xrightarrow{\ m\ } & G
\end{array}
\qquad \text{(associativity)} \qquad (2\mathrm{A}.2)
$$

$$
\begin{array}{ccc}
G \times T & \xrightarrow{\ \mathrm{id}\, \times e\ } & G \times G \\
\big\| {\scriptstyle =} & & \big\downarrow{\scriptstyle m} \\
G & = & G
\end{array}
\qquad
\begin{array}{ccc}
T \times G & \xrightarrow{\ e\, \times \mathrm{id}\ } & G \times G \\
\big\| {\scriptstyle =} & & \big\downarrow{\scriptstyle m} \\
G & = & G
\end{array}
\qquad \text{(unit)} \qquad (2\mathrm{A}.3)
$$

$$
\begin{array}{ccc}
G & \xrightarrow{\ (\mathrm{id},\, \iota)\ } & G \times G \\
\big\downarrow & & \big\downarrow{\scriptstyle m} \\
T & \xrightarrow{\ e\ } & G
\end{array}
\qquad
\begin{array}{ccc}
G & \xrightarrow{\ (\iota,\, \mathrm{id})\ } & G \times G \\
\big\downarrow & & \big\downarrow{\scriptstyle m} \\
T & \xrightarrow{\ e\ } & G
\end{array}
\qquad \text{(inverse)} \qquad (2\mathrm{A}.4)
$$

where the vertical arrow on the left hand side of the two diagrams (2A.4) is the unique morphism in the category  $\mathcal{C}$  to the terminal object and the vertical isomorphisms on the left of (2A.3) are the canonical isomorphisms of an object with the product of that object with the terminal object.

In the case of the category of groups this means that a group object is a group (with composition law denoted $+$ (though it is not clear yet that it is commutative) with a second composition law, denoted $\ast$ that is distributive over the first composition law in the sense that the following identity holds

$$(a + b) \ast (a' + b') = (a \ast a') + (b \ast b') \tag{2A.5}$$

This comes from the requirement that $\ast$ must be a morphism in the category **Group**. Let $0$ be the unit element for the composition law $+$ and $1$ the unit element for the composition law $\ast$. Putting $b = a' = 0$ in (2A.5) gives

$$a \ast b' = (a \ast 0) + (0 \ast b') \tag{2A.6}$$

On the other hand putting in $a' = b = 1$ in (2A.5) gives

$$a + b = (a + b) \ast (1 + 1)$$

and multiplying this with the inverse of $a + b$ for the star composition gives $1 = 1 + 1$ and hence $1 = 0$. Put this in (2A.6) to find that $a \ast b = a + b$ showing that the compositions are the same and then (2A.5) immediately gives that both are Abelian.

Thus a group object in the category of groups is Abelian and the second composition law is the same as the first.

Actually this can be proved more generally for monoid objects in the category of groups, [99 Kan].

There is a nice illustration of this in homotopy theory (and that is where the idea came from). This goes as follows. The second homotopy group, $\pi_2(X, \ast)$, of a based space $(X, \ast)$ is, as a set, the set of all homotopy classes of maps from the disk into $X$ that take the boundary circle into the base point $\ast$ of $X$.

For illustrational (and conceptual) purposes it is easier to think of homotopy classes of maps from the unit filled square to $X$ that take the boundary to the base point. Homotopically, of course, this makes no difference.

Now let

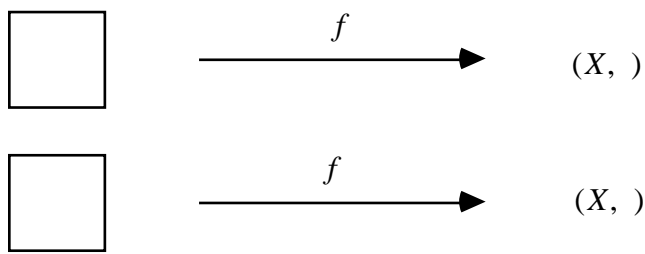

be two such maps. They can be glued together horizontally to give a map of the same kind (up to homotopy):

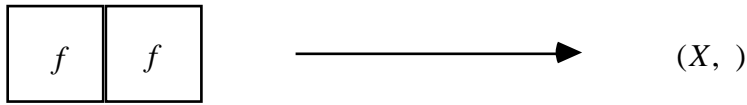

and this induces a composition on $\pi_2(X,\ )$ turning it into a group. Of course the two maps can also be glued together vertically, inducing another, a priori different, group structure.

Now take four such maps $f,\ f,\ g,\ g$. Then first gluing $f,\ f$ and $g,\ g$ together horizontally and then gluing the two results together vertically gives a map that can be depicted

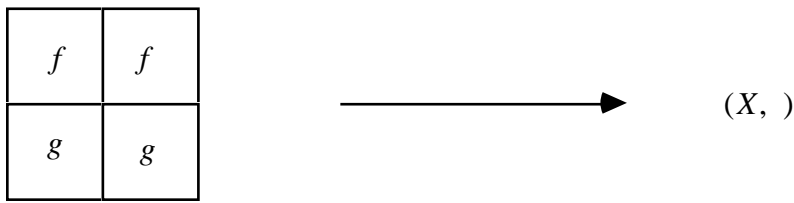

Obviously the same result is obtained by first gluing $f$ and $g$ together vertically, gluing $f$ and $g$ together vertically, and gluing the results together horizontally. This establishes the relation (2A,5) in the present case and shows that $\pi_2(X,\ )$ is Abelian. [3]

2A.7. *Comonoids in the category of groups.* Dually there is the notion of a cogroup object in a category.

For this let $\mathcal{C}$ be a category with direct sums and an initial object. Again the category of groups is an example with the one element group as initial object. The categorical direct sum in **Group** is what in group theory is called the free product.

A cogroup object in such a category $\mathcal{C}$ is an object $C\ \mathcal{C}$ together with a comultiplication $\mu: C\qquad C\ \ C$, a coinverse $\iota: C\qquad C$, and a counit morphism $\varepsilon: C\qquad I$. Here $I\ \mathcal{C}$ is the initial object and $\qquad$ stands for the direct sum in $\mathcal{C}$. These bits of structure are supposed to satisfy the dual axioms to those for a group object depicted by diagrams (2A.2) - (2A.4). that is the diagrams obtained by reversing all arrows (and replacing $m$ by $\mu$ and $e$ by $\varepsilon$) must be commutative. For a comonoid object leave out the coinverse and (the dual of) diagram (2A.4).

---

[3] In the present case of homotopy groups it can of course also easily be shown directy that vertical gluing and horizontal gluing give the same result and this is how things are done traditionallly in text books; see e..g. [94 Hu].

It is now a theorem, [99 Kan], that the underlying group of a comonoid object in the category of groups is free as a group.

This has much to do with the fact that the categorical direct sum in **Group** is given by the free product construction.

   2A.8. *Hopf's theorem on the cohomology of H-spaces*. An $H$-space is a based topological space $(X, )$ together with a continuous map $m: X \times X \qquad X$ such that $x \mapsto m(x, )$ and $x \mapsto m( , x)$ are homotopic to the identity. [4] The result of Heinz Hopf, [92 Hopf], see also [60 Félix et al.], alluded to is now as follows.

   Let $k$ be a field of characteristic zero and $X$ a path connected $H$-space such that $H (X;k)$ is of finite type then $H (X;k)$ is a free graded-commutative graded algebra.

Here 'finite type' means that each $H_i(X; k)$ is finite dimensional and the cohomology algebra is graded-commutative ( = commutative in the graded sense), i.e. $xy = (-1)^{\text{degree}(x)\text{degree}(y)} yx$. Thus the seemingly weak extra bit of structure '$H$-space' has a profound influence on the (cohomological) structure of a space.

   2A.9. *Intermezzo*: *Hopf algebras*. Let $R$ be a unital commutative ring. A graded module over $R$ is simply a collection of modules over $R$ indexed by the nonnegative integers [5]. Or, equivalently, it is a direct sum

$$M = \bigoplus_{i \in \mathbf{N} \{0\}} M_i \qquad\qquad (2A.10)$$

An element $x \quad M_i$ is said to be homogeneous of degree $i$. A graded module (2A.10) is said to be of finite type if each of the $M_i$ is of finite rank over the base ring $R$.

   The tensor product of two graded modules $M, N$ is graded by assigning degree $i + j$ to the elements from $M_i \quad N_j$.

A graded algebra over $R$ is a graded module (2A.10) equipped with a graded associative multiplication and a unit element

---

   [4] Often in the literature for an $H$-space it is also required that the 'multiplication' $m$ is associative up to homotopy. For the present result that is not required.

   [5] These will be the only kind of gradings occurring

$$m : M \otimes M \longrightarrow M, \quad m(M_i \otimes M_j) \subset M_{i+j}; \quad 1 \in M_0 \qquad (2A.11)$$

There are two notions of commutativity for graded algebras: (ordinary) commutativity, which means $xy = yx$, and graded-commutativity, which means $xy = (-1)^{\deg(x)\deg(y)} yx$. Both occur frequently in the literature and both will occur in the present paper [6]. 

    Correspondingly there are two versions for the multiplication in the tensor product of (the underlying graded modules) of graded rings, viz.

$$(x \otimes y)(x' \otimes y') = xx' \otimes yy'$$

$$\qquad (2A.12)$$

$$(x \otimes y)(x' \otimes y') = (-1)^{\deg(y)\deg(x')} xx' \otimes yy'$$

where in the second equation the elements $x, x', y, y'$ are supposed to be homogeneous. The sign factor in the second equation of (2A.12) is needed to ensure that the tensor product of two graded-commutative graded algebras is graded-commutative (as of course one wants it to be).

Dually a graded coalgebra over $R$ is a graded module equipped with a coassociative comultiplication and a counit

$$\mu : M \longrightarrow M \otimes M, \quad \mu(M_n) \subset \sum_{i+j=n} M_i \otimes M_j; \quad \varepsilon : M \longrightarrow R, \quad \varepsilon(M_i) = 0 \text{ for } i > 0$$

Just as in the algebra case there are two notions of cocommutativity and two ways to define a coalgebra structure on the tensor product of two graded coalgebras. These two are as follows. Let $C$ and $D$ be two graded coalgebras with comultiplications $\mu_C$, $\mu_D$. Write

$$\mu_C(x) = \sum x_i' \otimes x_i'', \quad \mu_D(y) = \sum y_j' \otimes y_j''$$

as sums of tensor products of homogeneous elements. Then the two graded coalgebra structures alluded to are

$$x \otimes y \mapsto \sum x_i' \otimes y_j' \otimes x_i'' \otimes y_j''$$

$$\qquad (2A.13)$$

$$x \otimes y \mapsto \sum (-1)^{\deg(y_j')\deg(x_i'')} x_i' \otimes y_j' \otimes x_i'' \otimes y_j''$$

---

[6] If all the odd degree summands of the graded ring are zero the two notions aggree. This can be used to unify things.

Next, a graded bialgebra $B$ is a comonoid object in the category of graded algebras or, equivalently, a monoid object in the category of graded coalgebras. Here again there are two versions depending on what algebra and coalgebra structures are taken on $B \otimes B$. First there is an 'ordinary' bialgebra which happens to carry a grading. In this case the algebra and coalgebra structures are given by the first formulas of (2A.12) and (2A.13). Second there is the 'grade-twist' version in which the algebra and coalgebra structures on the tensor product are given be the second formulas from (2A.12) and (2A.13). Here 'ordinary twist' and 'grade twist' respectively refer to the morphisms

$$x \otimes y \mapsto y \otimes x, \quad x \otimes y \mapsto (-1)^{\deg(x)\deg(y)} y \otimes x$$

which make their appearance when the conditions are written out explicitly in terms of diagrams or elements.

Finally, a graded Hopf algebra is a graded bialgebra that in addition carries a so-called antipode. That is a morphism $\iota$ of graded modules of degree 0 (so that $\iota(H_i) \subset H_i$) that satisfies

$$m(\mathrm{id} \otimes \iota)\mu = e\varepsilon \quad \text{and} \quad m(\iota \otimes \mathrm{id})\mu = e\varepsilon$$

A graded Hopf algebra over $R$ is connected if the grade zero part $H_0$ is equal to $R$ so that $e$ and $\varepsilon$ induce isomorphism of $R$ with $H_0$.

An element $x$ in a graded Hopf algebra (or bialgebra) is called primitive if it satisfies

$$\mu(x) = 1 \otimes x + x \otimes 1 \tag{2A.14}$$

These form a graded submodule $P(H)$ of the Hopf algebra $H$.

In the case of an 'ordinary twist' Hopf algebra the commutator product

$$[x, y] = xy - yx \tag{2A.15}$$

turns $P(H)$ into a Lie algebra (that happens to carry a grading such that the Lie bracket is of degree 0.

In the case of a 'graded twist' Hopf algebra take

$$[x, y] = xy - (-1)^{\deg(x)\deg(y)} yx \tag{2A.16}$$

to obtain a graded Lie algebra. That is a module equipped with a bilinear product $[ , ]$ that satisfies graded anticommutativity and the graded Jacobi identity:

$$[x,y] = (-1)^{\text{degree}(x)\text{degree}(y)}[y,x]$$

$$\text{(2A.17)}$$

$$[x,[y,z]] = [[x,y],z] + (-1)^{\text{degree}(x)\text{degree}(y)}[y,[x,z]]$$

2A.18. *Milnor-Moore theorem* (*topological incarnation*). Let $PX$ be the Moore path space of a path connected based topological space $(X, )$. That is the space of paths starting from   with specified length (which is what the adjective 'Moore' means in this context). Assigning to a path its endpoint defines a  continuous map $PX \quad X$, which is a fibration with   $X$, the space of Moore loops, as its fibre (over  ). As $PX$ is contractible the long exact homotopy sequence attached to this fibration gives isomorphisms $\pi_n(X) \qquad \pi_{n-1}( \ X)$. This can be used to transfer the Whitehead products $\pi_m(X) \times \pi_n(X) \qquad \pi_{m+n=1}(X)$  to a Lie product (of degree zero) $(\pi ( \ X) \quad k) \times (\pi ( \ X) \quad k) \quad ^{m \ x} \ \pi ( \ X) \quad k$, defining a graded Lie algebra $L_X$.

Composition of loops turns   $X$ into a topological monoid and, up to homotopy there is an inverse as well. Using the Alexander - Whitney and Eilenberg - Zilber chain complex equivalences, see [60 Félix, et al.], p. 53ff, and the fact that taking homology of chain complexes commutes with tensor products, ibid. p. 48, the composition $X \times \ X \qquad X$ and diagonal  $: X \qquad X \times \ X$ induce an algebra and coalgebra structure on $H ( \ X)$. Moreover, essentially because a loop in a product $X \times Y$ is a pair of loops and composition of loops seen this way goes component-wise, the comultiplication morphism $H ( \ X) \qquad H ( \ X) \quad H ( \ X)$ is an algebra morphism [7], ibid. p. 225.

All in all this turns   $H ( \ X)$ into a graded connected Hopf algebra (of the 'graded twist' kind).

Now let the coefficients ring used when taking cohomology be a field of characteristic zero.

2A.19. *Theorem* ([137 Milnor et al.], see also [60 Félix et al.], p. 293). Let $X$ be a simply connected path connected topological space. Then the Hurewicz homomorphism for   $X$ is an isomorphism of graded Lie algebras of   $L_X$ onto the graded Lie algebra of primitives of   $H ( \ X; k)$  and this isomorphism extends to an isomorphism of graded Hopf algebras of the universal enveloping algebra $UL_X$ with $H ( \ X; k)$. [8]

There is also a purely algebraic theorem that goes by the name 'Milnor-Moore theorem'. That one involves the notion of the universal enveloping algebra of a Lie algebra and will be discussed in subsection 2C below.

[7] This is the origin of the unfortunate but frequently used notation '  ' for the comultiplication in a Hopf algebra.

[8] Universal enveloping algebras are the topic of section 2B.1 below.

To conclude this section 2A let me briefly mention two more simple results that, I feel, qualify as 'niceness theorems'. Both say that the presence of a Hopf algebra (bialgebra) structure has implications for the underlying algebra.

2A.20. *Cartier's theorem on nilpotents in group schemes*. Let $H$ be a finite dimensional Hopf algebra over a field of characteristic zero. Then the underlying algebra has no nilpotents. Actually a much stronger statement holds, see [53 Dolgachev]. The usual statement is: A group scheme of finite type over a field of characteristic zero is smooth. See loc. cit. and [170 Voskresenskii], p. 7.

In characteristic $p > 0$, Cartier's theorem does not hold. On $k[X]/(X^p)$ where $k$ is a field of characteristic $p > 0$, there are the two comultiplications

$$X \mapsto 1 \quad X + X \quad 1, \quad X \mapsto 1 \quad X + X \quad 1 + X \quad X$$

and both define a bialgebra, and in fact Hopf algebra structure on $k[X]/(X^p)$. These two Hopf algebras (finite group schemes) are traditionally denoted $\alpha_p$ and $\mu_p$.

2A.21. Let $k$ be a field and $n$ an integer $\quad 2$. Then there is no bialgebra structure on the algebra $M^{n \times n}(k)$ of $n \times n$ matrices over $k$. See [49 Dascalescu et al.], p. 173.

It is a completely unknown which products of matrix algebras do carry (admit) a bialgebra structure.

Much of mathematics concerns statements as to what consequences follow from what assumptions. So it can be argued that there is nothing particularly special about the results described above. However, I feel that there is something special, something particularly elegant, about the results described. Part of the general problem is to understand why and in what sense.

Several of the theorems above are 'freeness theorems'. They say that in the presence of suitable extra structure an object is free. Here follow five more. For the first three the 'extra structure' is that the object in question is imbedded in a free object. In some categories that means nothing; in others it is a strong bit of extra structure. Just what categorical properties rule this behavior is completely unknown.

2A.22. *Nielsen–Schreier theorem*. A subgroup of a free group is free, [153 Schreier; 143 Nielsen]; [161 Suzuki], p. 181.

2A.23. *Shirshov–Witt theorem.* Lie subalgebras of a free Lie algebra are free, [158 Shirshov; 173 Witt]. There is also, up to a point, a braided version, [103 Kharchenko].

2A.24. *Bergman centralizer theorem.* The centralizer of a non–scalar element in a free power series ring $k\langle\!\langle X\rangle\!\rangle$ is of the form $k[[c]]$, [27 Bergman]; [47 Cohn], p. 244.. Here $c$ is a single element!

2A.25. The fundamental group of a cogroup object in the homotopy category of 'nice' based topological spaces is free. See [28 Berstein]. These objects are sometimes called $H-$spaces (as a kind of dual or opposite object to $H-$spaces).

2A.26. *Bott–Samelson theorem.* The homology algebra $H(\quad X;k)$ is a free algebra generated by $H(X;k)$, [34 Bott et al.; 28 Berstein].

Here is the suspension functor and is the loop space functor on based topological spaces. These are adjoint and there results a topological morphism $X \qquad X$. The multiplication comes from the fact that loops at the base point can be composed making a loop space an $H-$space.

## 2B. Universal object examples

Here the theme is that objects that are defined in terms of some universal property have a tendency to pick up extra bits a structure.

2B.1. *The universal enveloping algebra of a Lie algebra.* Let $A$ be a unital associative algebra over a unital commutative base algebra $R$. Associated to $A$ there is a Lie algebra structure on $A$ defined by the commutator difference

$$[x,y]_A = xy - yx \tag{2B.2}$$

Let $\mathfrak{g}$ be a Lie algebra. A Lie morphism from $\mathfrak{g}$ to a unital associative algebra $A$ is a module morphism $\varphi:\mathfrak{g} \qquad A$ such that $\varphi([x,y]_{\mathfrak{g}}) = [\varphi x, \varphi y]_A$. The universal enveloping algebra on $\mathfrak{g}$ is a unital associative algebra $U\mathfrak{g}$ together with a Lie morphism $i:\mathfrak{g} \qquad U\mathfrak{g}$ such that for each Lie morphism $\varphi:\mathfrak{g} \qquad A$ there is a unique morphism of associative algebras $\tilde{\varphi}:U\mathfrak{g} \qquad A$ such that $\tilde{\varphi} \circ i = \varphi$. Pictorially (in diagram form) this can be rendered as follows

$$\mathfrak{g} \longrightarrow U\mathfrak{g}$$

with $\varphi$ and $^1\tilde{\varphi}$ mapping to $A$

(2B.3)

The associative unital algebra $U\mathfrak{g}$ is a very nice one. For instance there is the Poincaré - Birkhoff - Witt theorem that specifies (under suitable circumstances) a monomial basis for it. This results basically from the construction of $U\mathfrak{g}$. (And one wonders whether this PBW theorem can be deduced directly from the characterizing universality property.)

What is of interest in the present setting is that the universality property immediately implies that $U\mathfrak{g}$ has more structure; in fact that it is a Hopf algebra. This arises as follows. Consider the associative algebra $U\mathfrak{g} \otimes U\mathfrak{g}$ and the morphism $x \mapsto 1 \otimes x + x \otimes 1$ from $\mathfrak{g}$ into it. It is immediate that this is a Lie morphism and hence there is a corresponding (unique) morphism of associative algebras $U\mathfrak{g} \to U\mathfrak{g} \otimes U\mathfrak{g}$. It is immediate that this turns $U\mathfrak{g}$ into a Hopf algebra.

There is a completely analogous picture for graded Lie algebras.

Of course the universal problem described here is an instance of an adjoint functor situation. Let **Lie** be the category of `Lie algebras (over $R$) and **Alg** the category of unital associative algebras (over $R$). Then associating to an associative algebra $A$ its commutator difference product is a (forgetful) functor $V : \mathbf{Alg} \to \mathbf{Lie}$ and $\mathfrak{g} \mapsto U\mathfrak{g}$ is a functor the other way that is left adjoint to it:

$$\mathbf{Lie}(\mathfrak{g}, V(A)) \cong \mathbf{Alg}(U\mathfrak{g}, A) \qquad (2B.4)$$

In the case of a forgetful functor a left adjoint to it yields what are often called free objects (as in this case). Thus $U\mathfrak{g}$ is the free associative algebra on the Lie algebra $\mathfrak{g}$.

A right adjoint functor to a forgetful functor gives cofree objects. An example of a cofree construction will occur below.

The very important notion of adjointness is due to Daniel Kan, [98 Kan] and as Saunders Mac Lane says in the preface of [126 Mac Lane] "Adjoint functors arise everywhere".

If $(F, G)$ is an adjoint functor pair, i.e. e.g. $\mathcal{C}(FX, Y) \cong \mathcal{D}(X, GY)$ functorialy (loosely formulated), one expects niceness properties for both the $FX$ 's and the $GY$ 's. And indeed many niceness results fall into this scope with the proviso that often these objects pick up extra properties which are not implicit in the adjoint situation alone.

2B.5 . *The group algebra of a group*, Much the same picture holds for the group algebra of a group. Except much easier. Here the 'forgetful functor' assigns to an algebra $A$ its group $A^*$ of invertible elements. Recall that the group algebra $kG$ of a group is the free module over $k$ with basis $G$ and the multiplication determined on this basis by the group multiplication. The adjointness equation now is:

$$\mathbf{Group}(G, A^*) \quad \mathbf{Alg}_k(kG, A) \qquad\qquad (2B.6)$$

There is again a Hopf algebra structure for free. For this, to put things formally on the same footing as in the case of the universal enveloping algebra, consider the morphism

$$G \quad (kG \quad kG)^* , \quad g \mapsto g \quad g$$

which by the adjointness equation (2B.6), gives rise to a morphism of algebras $kG \quad kG \quad kG$ turning $kG$ into a bialgebra (and a Hopf algebra using the group inverse). Of course in this case things are so simple that it is not worthwhile to go through this yoga.

2B.7. *Free algebras.* Everyone knows how to construct the free algebra over a module (or a set). The tensor algebra does the job and that is a very nice structure. Less known is that this also works in the setting **CoAlg** - **HopfAlg**, where **CoAlg** and **HopfAlg** are suitable categories of coalgebras and Hopf algebras over a suitable base ring. See [140 Moore] and [89 Hazewinkel]. This gives the free Hopf algebra on a coalgebra.

2B.8. *Cofree coalgebras.* Given a module $M$, the cofree coalgebra over [9] $M$ would be a coalgebra $C(M)$ together with a module morphism $C(M) \xrightarrow{\eta} M$ such that for each coalgebra $C$ together with a morphism of modules $C \xrightarrow{\alpha} M$ there is a unique morphism of coalgebras $\hat{\alpha}: C \quad C(M)$ such that $\eta\hat{\alpha} = \alpha$ .

---

[9] It pays to be terminologically careful in this contexr. I prefer to speak of the free algebra *on* a module and the cofree coalgebra *over* a module.

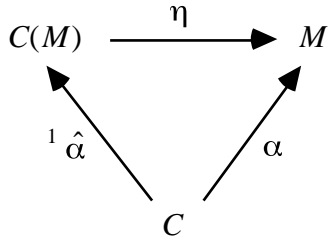

Whether the cofree coalgebra over a module always exist is not quite settled, [81 Hazewinkel]; they certainly exist in many cases. In the connected graded context they always exist and are given by the tensor coalgebra, again a very nice structure.

And in this connected graded context there is the **Alg - HopfAlf** version of the cofree Hopf algebra over an algebra, [89 Hazewinkel; 140 Moore].

**2B.9.** *The classifying spaces* $\mathbf{BU}_n$. A completely different kind of universal object is formed by the complex Grassmannians and their inductive limits the classifying spaces $\mathbf{BU}_n$.

Consider the complex vector space $\mathbf{C}^{n+r}$ and define the complex Grassmannian

$$\mathbf{Gr}_n(\mathbf{C}^{n+r}) = \{V: \ V \text{ is an } n \text{-dimensional subspace of } \mathbf{C}^{n+r}\} \qquad (2\text{B}.10)$$

This set has a natural structure of a smooth manifold (in fact a complex analytic manifold). Letting $r$ go to infinity (which technically means taking an inductive limit) gives the classifying space

$$\mathbf{BU}_n = \underrightarrow{\lim} \ \mathbf{Gr}_n(\mathbf{C}^{n+r}) = \mathbf{Gr}_n(\mathbf{C}^{\infty}) \qquad (2\text{B}.11)$$

It is also perfectly possible to define and work directly with the most right hand side of (2B.11). There is a canonical complex vector bundle over $\mathbf{BU}_n$ which is colloquially defined by saying the fibre over $x \in \mathbf{BU}_n$ "is $x$". More precisely this canonical vector bundle $\gamma_n$ is

$$\gamma_n = \{(x,v): \ x \in \mathbf{BU}_n, \ v \in x\} \text{ with projection } (x,v) \mapsto x, \ \gamma_n \longrightarrow \mathbf{BU}_n$$
(2B.12)

There is now the following universality/classifying property. For every paracompact space $X$ with an $n$-dimensional complex vector bundle $\xi$ over it there is a map $f_\xi: X \longrightarrow \mathbf{BU}_n$ such that $\xi$ is isomorphic (as a vector bundle) to the pullback $f_\xi^*(\gamma_n)$. Moreover $f_\xi$ is unique up to homotopy.

The remarkable thing here is that the classifying spaces $\mathbf{BU}_n$ are so elegant and simple (as are the universal bundles over them). There are more nice properties. Jumping the gun a little – these spaces will return later – the cohomology of these spaces is particularly nice

$$H\;(\mathbf{BU}_n;\mathbf{Z}) = \mathbf{Z}[c_1, c_2, \cdots, c_n], \;\; \deg(c_r) = 2r \qquad\qquad (2\mathrm{B}.13)$$

All this can be found in [96 Husemoller; 138 Milnor et al.] (and many other books).

## 2C. Niceness theorems for Hopf algebras

The structure of a Hopf algebra is a heavy one. Indeed at one time they were thought to be so rare that each and every one deserves the most careful study, {Kaplansky, #62}. This is not anymore the case. Hopf algebras abound. Still the structure is not strong enough to produce good niceness theorems. However if one adds conditions like graded and connected some strong structure theorems emerge. These are e.g. the Leray and Milnor - Moore theorems which will both be described immediately below. In addition there is the Zelevinsky theorem, a structure theorem due to Grünenfeld, {Grünenfeld, #63} and much more, see e.g. {Masuoka, 2007 #64}. However, whether the various available classification theorems for Hopf algebras qualify as niceness theorems is debatable. I think mostly not.

2C.1. *The Leray theorem on commutative Hopf algebras.* Let $H$ be a commutative graded connected Hopf algebra of finite type over a field of characteristic zero. Then the underlying algebra is commutative free. There is also a graded commutative version. The original theorem appears in [119 Leray]. For an up-to-date short account see [146 Patras]. There are all kinds of generalizations, e.g. to an operadic setting, see [145 Patras; 123 Livernet; 66 Fressé]

2C.2. *The Milnor - Moore theorem on cocommutative Hopf algebras.* Let $H$ be a cocommutative graded connected Hopf algebra of finite type over a field of characteristic zero. Then the underlying algebra is the universal enveloping algebra of the Lie algebra of primitives $P(H)$ of $H$, [137 Milnor et al.].

This is the algebraic incarnation referred to in 2A.18 above.

The Milnor-Moore theorem is a dual of the Leray theorem. To realize this recall from subsection 2B above that $U\mathfrak{g}$ is the free object in **Ass** on the object $\mathfrak{g}$ **Lie**.

## 2D. Large vs nice.

There is a tendency for (really) nice objects to be big (or very small). A prime

example is

$$\mathbf{Symm} = \mathbf{Z}[h_1, h_2, h_3, \cdots] \tag{2D.1}$$

the ring of polynomials over the integers in countably infinite many commuting variables over the integers. This object will be discussed in some detail further on.

Inversely big objects have a better change of being nice.

In this subsection I give some examples of this phenomenon.

2D.2. *Big projective models are free*. This result is due to Hyman Bass, [23 Bass]. For a precise statement see loc. cit. (corollary 3.2). The key ingredient is the following elegant observation [10].

If $P \quad Q \quad F$ with $F$ a non–finitely generated free module, then $P \quad F \quad F$.

The proof is simplicity itself and clearly shows the power and usefulness of infinity.

$$F \quad F \quad F \quad \cdots \quad P \quad Q \quad P \quad Q \quad \cdots$$
$$P \quad F \quad F \quad \cdots \quad P \quad F$$

2D.3. *General linear groups in various dimensions*. Let $k$ be the field of real numbers, complex numbers or even the quaternions. The general linear groups $\mathbf{GL}_n(k)$ for finite natural numbers are homotopically and cohomologically far from trivial.

Things change drastically in infinite dimension.

2D.4. *Kuiper's theorem*, [110 Kuiper]. Let $H$ be real or complex or quaternionic Hilbert space. Then the general linear group $\mathbf{GL}(H)$ is contractible.

There is also an important equivariant extension due to Graeme Segal, [157 Segal].

Much related is Bessaga's theorem, [29 Bessaga; 30 Bessaga et al.], to the effect that every infinite dimensional Hilbert space is diffeomorphic with its unit sphere.

Kuiper's famous theorem is the key to the classification of Hilbert manifolds, [37 Burghelea et al.; 56 Eells et al.; 57 Eells et al.; 141 Moulis; 142 Moulis].

2D.5. Here is a table on differential topology in various dimensions as things seem to be constituted at present.

---

[10] Hyman Bass calls it "an elegant little swindle".

| 1 | 2 | 3 | 4 | 5 | 6 | ··· ··· | < infinity | infinity |
|---|---|---|---|---|---|---|---|---|
| Real easy | Easy | difficult | difficult; boapw | good techniques | good techniques | | good techniques | real nice |

Here 'good techniques' refers mainly to Smale's handlebody theory. The acronym 'boapw' means 'best of all possible worlds' and refers to the fact that all $\mathbf{R}^n$ for $n$  4 have a unique differentiable structure, but $\mathbf{R}^4$ has over countably infinite different differentiable structures. [11]

### 2E. Extremal objects and niceness.

In the world of optimization theory and variational calculus and analysis it is relatively well known that extremal objects tend to be nice (have lots of symmetry), even when bifurcation occurs.

There are also various notions of minimality in algebra and topology and these also tend to be 'nice'. For instance the Sullivan minimal models for rational homotopy, see [60 Félix et al.], are definitely nice.

In the world of operads and PROP's etc. there are by way of example the following theorems, see [136 Merkulov].

– The minimal resolution of $\mathcal{Ass}$ is a differential graded *free* operad.

– The minimal resolution of $\mathcal{LieB}$ is a *free* differential graded PROP.

Sullivan minimal models and operads, PROP's etc are highly technical notions and giving details would take me far beyond the scope and intentions of this paper.

I have no doubt that there are more niceness results for minimal resolutions. [12]

### 2F. Uniqueness and rigidity and niceness

For instance **Symm**, see (2D.1) above and below, is unique and rigid as a coring object in the category of unital commutative rings and **MPR**, the Reutenauer-Malvenuto-Poirier Hopf algebra of permutations is rigid and likely unique, see [86 Hazewinkel; 88 Hazewinkel]. And indeed they are very nice objects.

---

[11] This is a fact that tends to make 'multiple world' enthousiasts happy.

[12] There are at least three meanings for the word 'resolution' and the phrase 'minimal resolution' in mathematics: resolution of singularities, resolution of a module in homological algebra, resolution in (automatic) theorem proving. Outside mathematics there are many more additional meanings.

## 2G. Counterexamples. and paradoxical objects

Not only objects and constructions can exhibit the 'niceness phenomenon' but also counterexamples. This subsection contains a few examples of that.

2G.1. *The Alexander horned sphere*. First the construction as illustrated by the picture below. Take a hollow cylinder closed at both ends and bend around so that the two ends face each other. Now from each end extrude a horn and interlock them as shown; there result two locations of disks facing each other. Repeat ad infinitum.

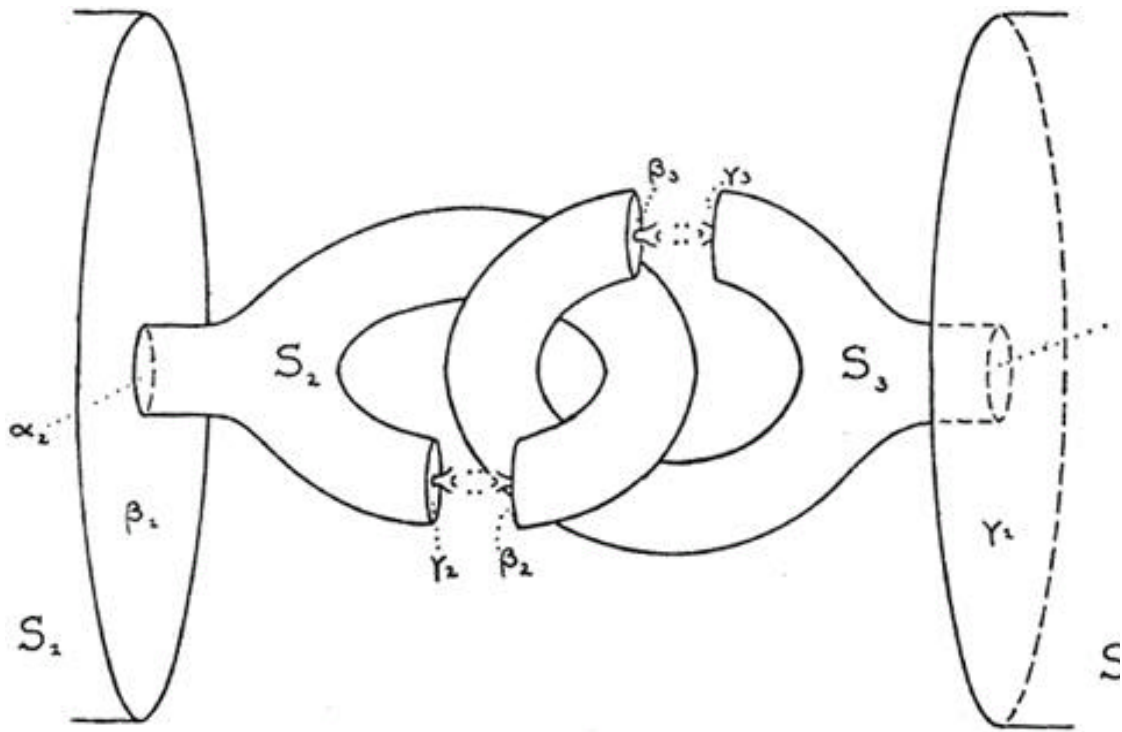

Fig. 1.  Alexander horned sphere

The Alexander horned sphere together with its interior is (homeomorphic to) a topological 3-ball. The exterior is not simply connected. This shows that the analogue of the Jordan-Schönflies theorem from dimension 2 does not hold in dimension 3.

For some more information on the Alexander horned sphere and its uses see [1; 2].

Somewhat surprisingly (to me in any case), the filled Alexander horned sphere can be used for a monohedral tiling of $\mathbf{R}^3$, [165 Tang].

2G.2. *The approximation property*. A Banach space is said to have the approximation property if every compact operator is a limit of finite rank operators.

Equivalently a Banach space $X$ has the approximation property if for every compact subset $K \quad X$ and every $\varepsilon > 0$ there is an operator $T: X \qquad X$ of finite rank such that $\|Tx - x\| < \varepsilon$ for all $x \quad K$.

Every Banach space with a (Schauder) basis has the approximation property. This includes Hilbert spaces and the $l^p$ spaces.

However, not every Banach space has the approximation property. In 1973 Per Enflo, [58 Enflo], constructed a counterexample.

I do not think this counterexample qualifies as a nice one. However the very nice Banach space of bounded operators on $l_2$ is also a counterexample, [162 Szankowski] [13].

2G.3. *The Banach-Tarski paradox*. In 1924 Stefan Banach and Alfred Tarski proved the following bizarre seeming statement, [22 Banach et al.].

For two bounded subsets $A$, $B$ of a Euclidean space of dimension at least three with nonempty interior there exist finite decompositions into disjoint subsets

$$A = A_1 \quad \cdots \quad A_k \qquad B = B_1 \quad \cdots \quad B_k$$

such that $A_i$ is congruent to $B_i$ for all $i = 1, \cdots, k$. I.e. $A_i$ becomes $B_i$ under a Euclidean motion.

This is now known as the strong form of the Banach-Tarski paradox. It does not hold in dimensions 1 and 2. A consequence is

A solid ball can be decomposed into a finite number of point sets that can be reassembled to form two balls identical to the original; see Fig. 2 below.

Here 'move' means a Euclidean space move: a combination of translations, rotations and reflections. For some more information on the Banach-Tarski paradox see [171 Wagon; 5; 6]

Thus 'move' is simple enough. The decomposition, however, is complicated. For one thing at least some of the components must be nonmeasurable. Also things are in three dimensions and Cantor-like sets in three dimensions are difficult to visualize. Fortunately Stan Wagon found a two dimensional analogue in hyperbolic space and the picture is remarkably beautiful; see Fig. 3 below.

---

[13] In the Wikipedia 'Szankowski' is rendered 'Shankovskii' which makes it quite hard to find the paper.

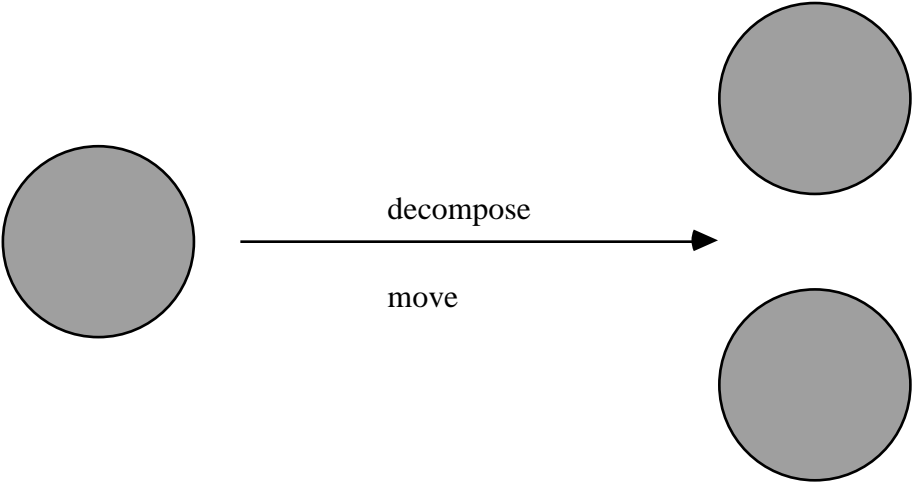

Fig. 2.  Banach-Tarski paradox: two balls out of one

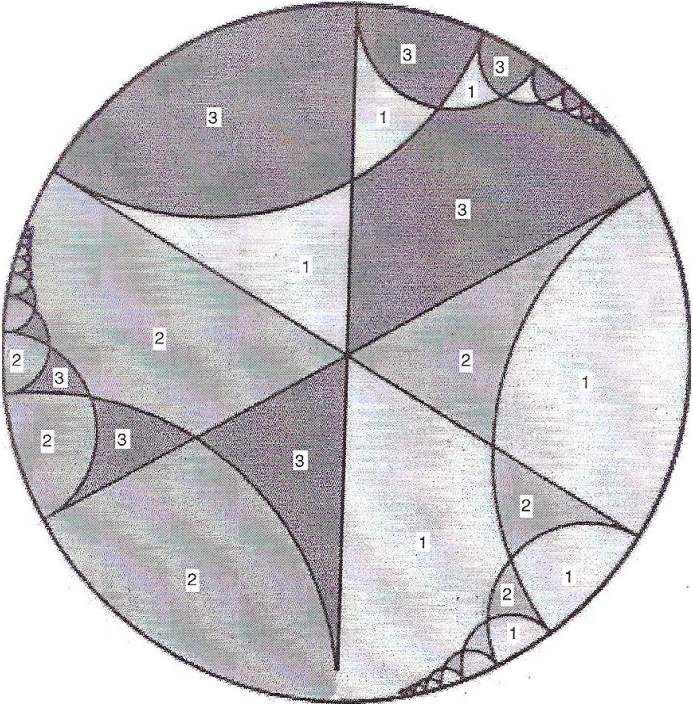

Fig. 3.  Banach-Tarski paradox: hyperbolic version

2G.4. *Julia and Fatou sets*. Here is a question not untypical of those that were asked in general (point-set) topology almost a century ago when people started to realize just how strange topological spaces could be.

Is it possible to divide the square into three regions so that the boundary between two of them is also the boundary between the two other pairs of regions.

The first answer was given by L E J Brouwer in the form of a simple construction repeated ad infinitum. However, the resulting picture is absolutely not beautiful. Nowadays there are the basins of attraction of discrete dynamical systems such as $x \mapsto x^3 - 1$ which has three basins of attraction (Fatou sets), one for each of the roots of $x^3 - 1$ and each pair has the same boundary (Julia set), see [11].

This is part of the world of fractals and (deterministic) chaos, [156 Schuster], and many of the pictures are extraordinarily beautiful, [147 Peitgen et al.] [14].

2G.5. *Sorgenfrey line*. As a set the Sorgenfrey is the set of real numbers. It is given a topology by taking as a basis the halfopen intervals $[a,b)$, $a < b$. This topology is finer than the usual one. For instance the sequence $\{n^{-1}\}_{n\ \mathbf{N}}$ converges to zero but $\{-n^{-1}\}_{n\ \mathbf{N}}$ does not. The Sorgenfrey line serves as a counterexample to several topological properties, [160 Steen et al.]. The point here (as far as this paper is concerned) is not that such counterexamples exist but that there is such a nice regular one. There is also a Sorgenfrey plane, loc. cit.

For some more information see also [125 Lukes; 13].

2G.6. *Exotic spheres*. A further example that fits in this section is that of exotic spheres (Milnor spheres). This deals with existence of differentiable structures on topological spheres, especially the seven dimensional ones, that differ from the standard one. They were the first examples of this phenomenon of distinct differentiable structures on the same topological manifold. This topic is rather more technical, and so I content myself with giving two references to internet accessible documents, [10; 152 Rudyak].

## 2H. An excursion into formal group theory.

A one dimensional formal group law over a commutative unital ring $A$ is a power series $F(X,Y)$ in two variables with coefficients in $A$ such that

$$F(X,0) = X, \ F(0,Y) = Y, \ F(X,F(Y,Z)) = F(F(X,Y),Z) \qquad (2H.1)$$

Two examples are the multiplicative formal group law and the additive formal group law

[14] Beautiful and arresting enough that the Sparkasse in Bremen organized an exhibition of them in 1984.

$$\hat{\mathbf{G}}_m(X,Y) = X + Y + XY, \quad \hat{\mathbf{G}}_a(X,Y) = X + Y \qquad (2\text{H}.2)$$

Both examples are nontypical in that they are polynomial; polynomial formal group laws are very rare.

More generally for any $n$, including $n = \infty$, an $n$–dimensional formal group over $A$ is an $n$–tuple of power series in two groups of $n$ indeterminates $F(X;Y)$ such that

$$F(X;0) = X, \ F(0;Y) = Y, \ F(X;F(Y;Z)) = F(F(X;Y);Z) \qquad (2\text{H}.3)$$

However, certainly from the point of view of applications, one dimensional formal groups are by far the most important, especially one dimensional formal groups over the integers , rings of integers of algebraic number fields, and over polynomial rings over the integers.

The only other that currently seems important is the infinite dimensional formal group $\hat{W}$ of the Witt vectors which is defined by the same polynomials that define the addition of Witt vectors; see the next subsection 2I.

A standard reference for formal groups is [80 Hazewinkel].

2H.4. *Lazard commutativity theorem*. Let $A$ be a ring that has no elements that are simultaneously torsion and nilpotent. Then every one dimensional formal group over $A$ is commutative; i.e. satisfies $F(X,Y) = F(Y,X)$.

2H.5. *Universal formal groups*. Given a formal group $F(X,Y)$ over $A$ and a morphism of rings $\alpha: A \longrightarrow B$ one obtains a formal group $\alpha F(X,Y)$ over $B$ by applying $\alpha$ to the coefficients of $F(X,Y)$.

A one dimensional commutative formal group $F_L(X,Y)$ over a ring $L$ is called universal [15] if for every one dimensional formal group $F(X,Y)$ over a ring $A$ there is a unique morphism of rings $\alpha^F: L \longrightarrow A$ such that $\alpha^F F_L(X,Y) = F(X,Y)$.

That such a thing exists and is unique is a triviality. What is very remarkable is the theorem of Lazard, [114 Lazard], that $L$ is the ring of polynomials in an infinity of indeterminates over the integers. The standard proof is a bitch and highly computational.

2H.6. *Morphisms*. A morphism of formal groups from an $m$–dimensional formal group $F(X;Y)$ to an $n$–dimensional formal group $G(X;Y)$ is an $n$–tuple of power series in $m$ indeterminates $\varphi(X)$ such that

$$\varphi(0) = 0, \quad G(\varphi(X); \varphi(Y)) = \varphi(F(X; Y))$$

---

[15] This is a rather different 'universal' than e.g. in 'universal enveloping algebra'. The 'L' in these sentences stands for Lazard.

If $\varphi(X) \equiv X \mod (\text{degree } 2)$ the morphism is said to be strict.

2H.7. *Logarithms*. Let $A$ be a ring of characteristic zero so that the canonical ring morphism $A \to A \otimes_{\mathbf{Z}} \mathbf{Q} = A_{\mathbf{Q}}$ is injective; let $F(X, Y)$ be a one dimensional formal group over $A$. Then over $A_{\mathbf{Q}}$ there exists a power series $f(X) = X + a_2 X^2 + \cdots$ such that

$$F(X, Y) = f^{-1}(f(X) + f(Y)) \qquad (2\text{H}.8)$$

Here $f^{-1}$ is the compositional inverse of $f$, i.e. $f^{-1}(f(X)) = X$. This $f$ is called the logarithm of $F$. In the case of the multiplicative formal group, see (2H.2), the logarithm is

$$\log(1 + X) = X - 2^{-1} X^2 + 3^{-1} X^3 - 4^{-1} X^4 + \cdots$$

Indeed, $\log(1 + X + Y + XY) = \log(1 + X) + \log(1 + Y)$. The terminology derives from this example. The logarithm of a formal group is a strict isomorphism of the formal group to the additive formal group; but over $A_{\mathbf{Q}}$.

It is at the level of logarithms that the recursive structure of formal groups appears; a recursive structure that was totally unexpected.

There are also logarithms for more dimensional commutative formal groups.

2H.9. *p-typical formal groups*. A one dimensional formal group over a characteristic 0 ring is *p*–typical if its logarithm is of the form

$$f(X) = X + b_1 X^p + b_2 X^{p^2} + \cdots$$

There is a better definition, see [80 Hazewinkel], which works always and also in the more dimensional case. But this one will do for the purposes of the present paper.

Over a $\mathbf{Z}_{(p)}$–algebra every formal group is strictly isomorphic to a *p*–typical one, [[39 Cartier]. If the ring over which the formal group is defined is of characteristic zero the isomorphism is easily described: take the logarithm and change all coefficients of non-*p*-powers of $X$ to zero.

2H.10. *The universal p–typical formal group*, [82 Hazewinkel]. Take a prime number $p$ and consider the following ring with endomorphism

$$\mathbf{Z}[V] = \mathbf{Z}[V_1, V_2, V_3, \cdots], \quad \psi(V_n) = V_n^p \tag{2H.11}$$

Define

$$a_n(V) = \sum_{i_1 + \cdots + i_r = n} p^{-r} V_{i_1} V_{i_2}^{p^{i_1}} V_{i_3}^{p^{i_1 + i_2}} \cdots V_{i_r}^{p^{i_1 + \cdots + i_{r-1}}} \tag{2H.12}$$

Thus the first few of these polynomials are

$$a_1(V) = p^{-1} V_1, \quad a_2(V) = p^{-2} V_1 V_1^p + p^{-1} V_2,$$

$$a_3(V) = p^{-3} V_1 V_1^p V_1^{p^2} + p^{-2} V_1 V_2^p + p^{-2} V_2 V_1^{p^2} + p^{-1} V_3$$

This sequence of polynomials has both a left and a right recursive structure.

The left recursive structure is

$$a_n(V) = \sum_{i=1}^{n} p^{-1} V_i \psi^i(a_{n-i}(V)) \quad \text{(where } a_0(V) = 1) $$

and the right recursive structure is

$$p a_n(V) = a_{n-1}(V) V_1^{p^{n-1}} + a_{n-2}(V) V_2^{p^{n-2}} + \cdots + a_1(V) V_{n-1}^p + V_n$$

Now consider

$$f_V(X) = X + a_1(V) X^p + a_2(V) X^{p^2} + a_3(V) X^{p^3} + \cdots$$

$$F_V(X, Y) = f_V^{-1}(f_V(X) + f_V(Y)) \tag{2H.13}$$

The left recursive structure is used to prove that $F_V(X, Y)$ is integral, i.e. has its coefficients in $\mathbf{Z}[V]$ and hence is a formal group over $\mathbf{Z}[V]$ and, subsequently, to prove that it is the universal $p$-typical formal group which means that every $p$-typical formal group can be obtained from it by a suitable ring morphism from $\mathbf{Z}[V]$.

The right recursive structure then leads to important applications to e.g. complex cobordism theory in algebraic topology and Dirichlet series in number theory.

The important thing here is not that a universal $p$–typical formal group exists but that it has these very simple and elegant recursive structures.

The universal $p$–typical formal groups can be simply fitted together to give a construction of the universal formal group.

2H.14. *Formal groups from cohomology*. Let $h^\cdot$ be a multiplicative extraordinary cohomology theory with first Chern classes. What all these words really mean is not so important at the present stage. Suffice that many of the better known cohomology theories are like this. The point is that under these circumstances there is a universal formula for the first Chern class of a tensor product of line bundles in terms of the first Chern classes of the factors.

$$c_1(\xi \otimes \eta) = \sum_{i,j} a_{ij} c_1(\xi)^i c_1(\eta)^j$$

defining a formal group over $h^\cdot(pt)$.

$$F_{h^\cdot}(X,Y) = \sum a_{ij} X^i Y^j, \quad a_{ij} \in h^\cdot(pt)$$

Here are some examples.

$h^\cdot = H^\cdot$, ordinary cohomology, $F_H = \hat{\mathbf{G}}_a$, the additive formal group

$h^\cdot = K^\cdot$, complex $K$–theory, $F_K(X,Y) = X + Y + uXY$, where $u$ is the Bott periodicity element; a version of the multiplicative formal group.

$h^\cdot = \mathbf{MU}^\cdot$, complex cobordism. In this case the formal group has logarithm $f_{\mathbf{MU}}(X) = \sum_{n=0}^{\infty} \frac{[\mathbf{CP}^n]}{n+1} X^{n+1}$. Here $\mathbf{CP}^n$ is $n$–dimensional complex projective space and $[\mathbf{CP}^n]$ is its complex cobordism class in $\mathbf{MU}^\cdot(pt)$. This profound result is due to A S Mishchenko, see appendix 1 of [144 Novikov].

$h^\cdot = \mathbf{BP}^\cdot$, Brown–Peterson cohomology, the 'prime $p$ part' of complex cobordism. Its formal group is the $p$–typification of the one of complex cobordism, so that its logarithm is $f_{\mathbf{BP}}(X) = \sum_{r=0}^{\infty} \frac{[\mathbf{CP}^{p^r-1}]}{p^r} X^{p^r}$.

For more details see [80 Hazewinkel] and the references given there and especially [151 Ravenel].

There is more. The formal group of complex cobordism is the universal one, [149 Quillen].

The remarkable, elegant and nice aspect here is that in terms of cobordism the universal formal group is so simple and regular.

It follows from the Quillen theorem that $F_{\mathbf{BP}}(X,Y)$ with logarithm $f_{\mathbf{BP}}(X)$ is the universal $p$–typical formal group law. But there is also an explicit construction of the universal $p$–typical formal group law, (2H.13). This has all kinds of consequences for complex cobordism and Brown–Peterson cohomology, see [84 Hazewinkel; 80 Hazewinkel; 151 Ravenel]

Quillen's theorem also goes a fair way towards establishing that complex cobordism is the most general cohomology theory.

### 2I. The amazing Witt vectors and their gracious applications [16]

Let **CRing** be the category of unital commutative associative rings. The big Witt vectors constitute a functor $W : \mathbf{CRing} \longrightarrow \mathbf{CRing}$ which has an amazing number of universality properties. For a fair amount of information on this functor see [89 Hazewinkel] and the references quoted there.

2I.1. *Definition of the functor of the big Witt vectors*. As a set $W(A) = (A)$ is the set of all power series with coefficients in $A$ with constant term 1.

$$W(A) = (A) = \{1 + a_1 t + a_2 t^2 + a_3 t^3 + \cdots : a_i \in A\} \tag{2I.2}$$

Multiplication of such power series defines an Abelian group structure on $W(A)$ with as neutral element the power series 1. This is the underlying group of the to be defined ring structure on $W(A)$. The multiplication on $W(A)$ is uniquely determined by the requirement that the very special power series $(1 - xt)^{-1}$ multiply as

$$(1 - xt)^{-1} \ast (1 - yt)^{-1} = (1 - xyt)^{-1} \tag{2I.3}$$

and the demands of distributivity (of multiplication over addition) and functoriality. Just how this works will be indicated immediately below.

The functoriality of $W(-)$ is component-wise, i.e. it is given by

$$W(f)(1 + a_1 t + a_2 t^2 + a_3 t^3 + \cdots) = 1 + f(a_1)t + f(a_2)t^2 + f(a_3)t^3 + \cdots \tag{2I.4}$$

The functor $W$ is obviously representable by the ring $\mathbf{Symm} = \mathbf{Z}[h_1, h_2, h_3, \cdots]$ of polynomials in a countable infinity of indeterminates over the integers. The functorial

correspondence is:

$$1 + a_1 t + a_2 t^2 + a_3 t^3 + \cdots \qquad f: Symm \longrightarrow A, \; f(h_n) = a_n \qquad (2I.5)$$

It is convenient to view the $h_n$ as the complete symmetric functions in another countably infinite set of indeterminates $\xi_1, \xi_2, \xi_3, \cdots$ which can be encoded as

$$1 + h_1 t + h_2 t^2 + h_3 t^3 + \cdots = \prod_i \frac{1}{(1 - \xi_i t)} \qquad (2I.6)$$

Now let $h_1, h_2, h_3, \cdots$ be a second set of commuting indeterminates viewed as the complete symmetric functions in $\eta_1, \eta_2, \eta_3, \cdots$ that commute with the $\xi$. Then distributivity requires that

$$(1 + h_1 t + h_2 t^2 + h_3 t^3 + \cdots)(1 + h_1 t + h_2 t^2 + h_3 t^3 + \cdots) = \prod_{i,j} \frac{1}{(1 - \xi_i \eta_j t)} \quad (2I.7)$$

This makes sense because the right hand side of (2I.7) is symmetric in the $\xi$ and in the $\eta$ and so, by the fundamental symmetric functions theorem there are unique polynomials

$$\wp_1(h_1; h_1), \quad \wp_2(h_1, h_2; h_1, h_2), \quad \wp_3(h_1, h_2, h_3; h_1, h_2, h_3), \; \cdots \qquad (2I.8)$$

such that

$$(1 + h_1 t + h_2 t^2 + h_3 t^3 + \cdots)(1 + h_1 t + h_2 t^2 + h_3 t^3 + \cdots)$$
$$= 1 + \wp_1(h_1; h_1)t + \wp_2(h_1, h_2; h_1, h_2)t^2 + \wp_3(h_1, h_2, h_3; h_1, h_2, h_3)t^3 + \cdots \qquad (2I.9)$$

(That the multiplication polynomials $\wp_n$ depend only on the first $n$ $h_i$ and $h_i$ is easily seen by degree considerations.)

By functoriality these polynomials determine the multiplication on each $W(A)$ in the sense that for $a(t) = 1 + a_1 t + a_2 t^2 + a_3 t^3 + \cdots$ and $b(t) = 1 + b_1 t + b_2 t^2 + b_3 t^3 + \cdots$ in $W(A)$ their product is

$$a(t) \cdot b(t) = 1 + \wp_1(a_1; b_1)t + \wp_2(a_1, a_2; b_1, b_2)t^2 + \wp_3(a_1, a_2, a_3; b_1, b_2, b_3)t^3 + \cdots$$

Of course the sum in $W(A)$ is also defined by universal polynomials. These are

$$\Sigma_n(h_1, \cdots, h_n; h_1, \cdots, h_n) = \sum_{i+j=n} h_i h_j \quad \text{where } h_0 = h_0 = 1 \qquad (2I.10)$$

Another way of expressing most of this is to say that

$$h_n \mapsto \,_n(h_1 \quad 1, \cdots, h_n \quad 1; 1 \quad h_1, \cdots, 1 \quad h_n)$$

$$h_n \mapsto \,_n(h_1 \quad 1, \cdots, h_n \quad 1; 1 \quad h_1, \cdots, 1 \quad h_n)$$

(2I.11)

define on **Symm** (most of) the structure of a coring object in the category **CRing**, which hence, via (2I.5) defines a functorial ring structure on the $W(A)$.

2I.12. *Lambda rings and sigma rings.* A pre-sigma-ring (pre-$\sigma$-ring) is a unital commutative ring $A$ that comes with extra nonlinear operators that behave (in a very real sense) like symmetric powers. That is, there are operators

$$\sigma_i : A \qquad A, \; i = 1, 2, \cdots \; ; \; \sigma_1 = \mathrm{id}$$

(2I.13)

such that

$$\sigma_n(x + y) = \sigma_n(x) + \sum_{i=1}^{n-1} \sigma_i(x) \sigma_{n-i}(y) + \sigma_n(y)$$

(2I.14)

It is often useful to have the notation $\sigma_0$ for the operator that takes the constant value 1. This notion is equivalent to the better known one of a pre-lambda-ring (pre-$\lambda$-ring) but works out just a bit better notationally. The two sets of operations are related by the Wronski-like relations

$$\sum_{i=0}^{n} (-1)^i \sigma_i(x) \lambda_{n-i}(x) \; = \; 0$$

The lambda operations behave like exterior powers.

Let $\varphi : A \qquad B$ be a morphism in **CRing** and let both $A$ and $B$ carry pre-sigma-ring structures. Then the morphism is said to be a morphism of pre-sigma-rings if it commutes with the sigma operations, i.e. $\varphi(\sigma_n^A(x) = \sigma_n^B(\varphi(x))$. A pre-sigma-ring is a sigma ring if the operations satisfy certain universal formulas when iterated and when applied to a product. This is conveniently formulated as follows.

Consider the ring of big Witt vectors $W(A)$ and write an element of it (formally) as

$$a(t) \; = \; 1 + a_1 t + a_2 t^2 + a_3 t^3 + \cdots \; = \; \prod_i \frac{1}{(1 - \xi_i t)}$$

Then

$$\sigma_n(a(t)) = \sum_{i_1 \quad i_2 \quad \cdots \quad i_n} (1 - \xi_{i_1}\xi_{i_2}\cdots\xi_{i_n}t)^{-1} \tag{2I.15}$$

(when written out in terms of the $a_i$ which can be done by the usual symmetric function yoga). This defines a pre-sigma-ring structure on $W(A)$

A pre-sigma-ring $A$ is a sigma-ring if

$$\sigma_t : A \longrightarrow W(A), \quad x \mapsto 1 + \sigma_1(x)t + \sigma_2(x)t^2 + \sigma_3(x)t^3 + \cdots$$

is a morphism of pre-sigma rings. It is a theorem that $W(A)$ is in fact a sigma-ring. This involves the study of the morphism

$$\sigma^{W(A)} : W(A) \longrightarrow W(W(A)) \tag{2I.16}$$

which I like to call the Artin-Hasse exponential [17].

A ring morphism between sigma-rings is a sigma-ring morphism if it is a morphism of pre-sigma-rings. Let **SigmaRing** be the category of sigma-rings.

Let

$$s_1 : W(A) \longrightarrow A, \quad a(t) \mapsto a_1 \tag{2I.17}$$

be the morphism of rings that assigns to a 1-power-series its first coefficient. The Witt vectors now have the following universality property. Let $S$ be a sigma-ring, $A$ a ring and $\varphi : S \longrightarrow A$ a morphism of rings, then there is a unique morphism of sigma-rings $\hat{\varphi} : S \longrightarrow W(A)$ such that the following diagram commutes

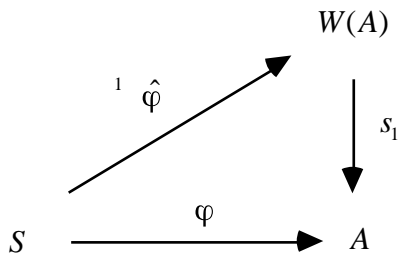

---

So $W(A) \to A$ is the cofree sigma-ring over the ring $A$. Or in other words the functor $W(-)$: **Cring** $\to$ **SigmaRing** is right adjoint to the functor the other way that forgets about the sigma structure.

2I.18. *The comonad structure on the big Witt vectors*. A comonad (also called cotriple) $(T, \mu, \varepsilon)$ in a category $\mathcal{C}$ is an endo functor $T$ of $\mathcal{C}$ together with a morphism of functors $\mu$: $T \to TT$ and a morphism of functors $\varepsilon$: $T \to$ id such that

$$(T\mu)\mu = (\mu T)\mu, \quad (\varepsilon T)\mu = \mathrm{id} = (T\varepsilon)\mu \qquad (2I.19)$$

And a coalgebra for the comonad $(T, \mu, \varepsilon)$ is an object in the category $\mathcal{C}$ together with a morphism $\sigma$: $C \to TC$ such that

$$\varepsilon_C \sigma = \mathrm{id}, \quad (T(\sigma))\sigma = (\mu_{TC})\sigma \qquad (2I.20)$$

It is now a theorem, [89 Hazewinkel] that the Artin -Hasse exponential (2I.16), which is functorial, together with the functorial morphism (2I.17) form a cotriple and that the coalgebras for this cotriple are precisely the sigma-rings.

2I.21. *The sigma and lambda ring structures on* **Symm**. Consider

$$\mathbf{Symm} = \mathbf{Z}[h_1, h_2, h_3, \cdots] \subset \mathbf{Z}[\xi_1, \xi_2, \xi_3, \cdots] \qquad (2I.22)$$

as before. There is a unique sigma-ring structure on $\mathbf{Z}[\xi]$ determined by

$$\sigma_n(\xi_i) = \xi_i^n \qquad (2I.23)$$

(The corresponding lambda operations are $\lambda_1(\xi_i) = \xi_i$, $\lambda_n(\xi_i) = 0$ for $n \geq 2$ so that the $\xi_i$ are like line bundles and this is a good way of thinking about them.) The subring **Symm** is stable under these operations and so there is an induced sigma-ring structure on **Symm**.

It is now a theorem that **Symm** with this particular sigma-ring structure is the free sigma-ring on one generator. More precisely:

For every sigma ring $S$ and element $x \in S$ there is a unique morphism of sigma-rings **Symm** $\to S$ that takes $h_1$ into $x$.

The universality properties described in 2I.12, 2I.18, 2I.21 are far from unrelated; see section 3C below.

A totally different universality property of the Witt vectors is the following one,

2I.24. *Cartier's first theorem*. The (infinite dimensional) formal group of the Witt vectors 'is' the sequence of addition polynomials $\Sigma_1, \Sigma_2, \cdots$ in $X_1, X_2, \cdots; Y_1, Y_2, \cdots$. This formal group is denoted $\hat{W}$. A fourth universality property off the Witt vectors holds in this setting.

Given two formal groups $F$ and $G$ of dimensions $m$ and $n$ respectively a morphism of formal groups $\alpha: F \longrightarrow G$ is an $n$-tuple of power series with zero constant terms $\alpha_1, \cdots \alpha_n$ in $m$ variables such that

$$G(\alpha_1(X), \cdots, \alpha_n(X); \alpha_1(Y), \cdots, \alpha_n(Y)) = (\alpha_1(F(X,Y)), \cdots, \alpha_n(F(X,Y)))    (2I.25)$$

A curve in an $n$-dimensional formal group $F$ is simply an $n$-tuple of power series in one variable, say, $t$. In $\hat{W}$ consider the particular curve $\gamma_0(t) = (t, 0, 0, \cdots)$. Then Cartier's first theorem says that for every formal group $F$ and curve $\gamma(t)$ in it there is a unique morphism of formal groups $\hat{W} \longrightarrow F$ that takes $\gamma_0(t)$ into $\gamma(t)$.

## 2J. The star example: Symm.

Here is a list of most of the objects with which this subsection will be concerned. Those which have not already been defined above will be described in section **3C** below.

— **Symm** $= \mathbf{Z}[h_1, h_2, \cdots] = \mathbf{Z}[c_1, c_2, \cdots] \subset \mathbf{Z}[\xi_1, \xi_2, \cdots]$, the ring of symmetric functions in an infinity of indeterminates. Here $h_n$ is the $n$-th complete symmetric function in the $\xi$ 's and the $c_n$ stand for the elementary symmetric functions. I am writing $c_n$ rather than $e_n$ because in the present context the $c_n$ will correspond to Chern classes.

— $U(\mathbf{Z})$, the universal lambda ring on one generator

— $R(W)$, The representing ring of the functor of the big Witt vectors; see subsection 2I above.

— $R(S) = \bigoplus\limits_{n=0}^{\infty} R(S_n)$, the direct sum of the rings of (the Grothendieck groups of) complex representations of the symmetric groups with the so-called exterior product; if $\rho$ is a representation of $S_r$ and $\sigma$ is a representation of the symmetric group on $s$ letters $S_s$ then $\rho\sigma = \mathrm{Ind}_{S_r \times S_s}^{S_{r+s}}(\rho \times \sigma)$. By decree ' $R(S_0)$' is equal to $\mathbf{Z}$. There is also a comultiplication: if $\sigma$ is a representation of $S_n$   $\mu(\sigma) = \operatorname*{Res}_{S_r \times S_s}^{S_n}(\sigma)$ . Together with $r+s=n$ obvious unit and counit morphisms this defines a Hopf algebra. (The antipode comes for free because of the graded connected situation.)

— $R_{\mathrm{rat}}(\mathrm{GL}\ )$ , the (Grothendieck) ring of rational representations of the infinite linear group.

— $E(\mathbf{Z})$ , the value of the exponential functor from [91 Hoffman] on the ring of integers.

— $U(\hat{W})$ , the covariant bialgebra of the formal group of the Witt vectors.

— $H\ (\mathbf{BU};\mathbf{Z})$ , the cohomology of the classifying space of complex vector bundles, **BU.**

— $H\ (\mathbf{BU};\mathbf{Z})$ , the homology of the classifying space **BU.**

These are all isomorphic and that implies that **Symm** is very rich in structure indeed. Nor is that all. For instance each of the components $R(S_n)$ of $R(S)$ is a lambda ring in its own right (inner plethysm). Further the functor of the big Witt vectors is lambda ring valued. However, this paper is not about **Symm** and its extraordinarily rich structure [18], but about niceness results. That includes 'nice proofs'. That is proofs of isomorphism between all these objects that derive from their universality, (co)freeness, ... properties and rely minimally on calculations. To what extent there are currently such proofs will be discussed below in subsection **3C**.

Two more objects that fit in this picture are the rational Witt vector functor in its role in the $K$-theory of endomorphisms, [16 Almkvist], and the Grothendieck group $K(\mathbf{P}_A)$ of polynomial functors $\mathbf{Mod}_A \qquad \mathbf{Mod}_k$ , where $A$ is an algebra over a field $k$ and $_A\mathbf{Mod}$ is the category of right $A$–modules, [127 Macdonald].. If $A = k$ this object is again isomorphic to the nine objects listed above.

The various isomorphisms and relations concerning whom I think I have something to say are depicted in the diagram below.

---

[18] I plan a future paper on that; meanwhile see [89 Hazewinkel].

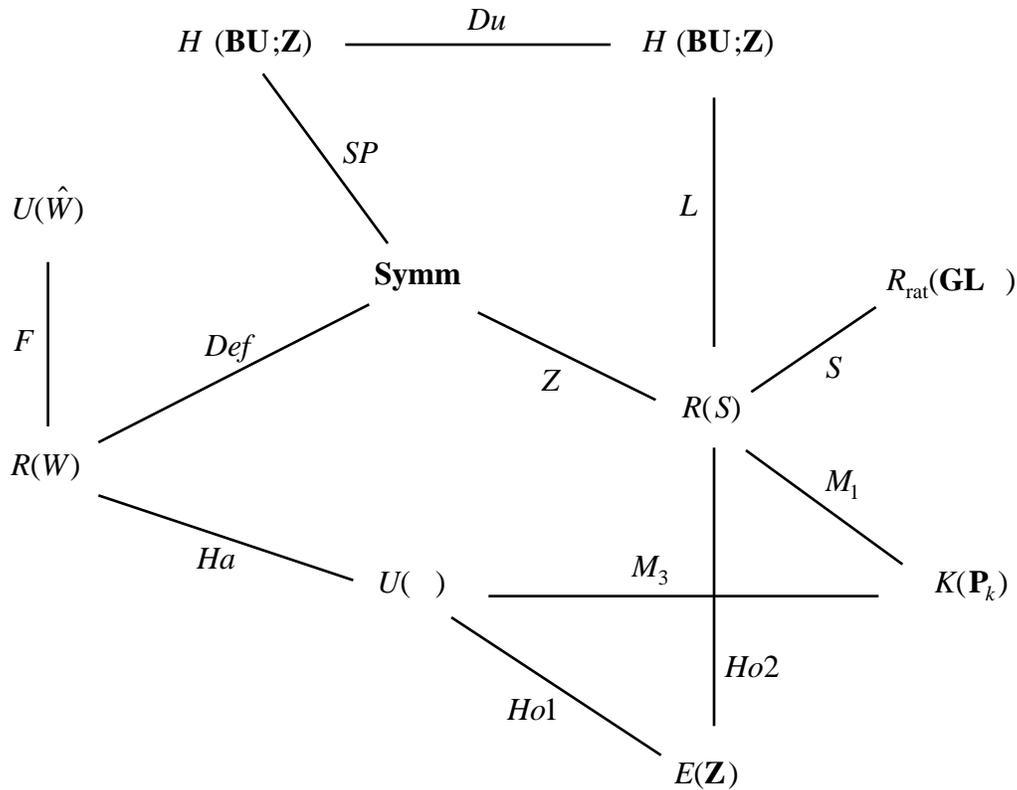

The bottom object here, viz $E(\mathbf{Z})$, has not yet been described in any way. It is again defined by an adjoint functor situation and, again, it is one which picks up extra structure. It will be described and discussed briefly in section **3C** below.

Also it seems from the diagram that the Hopf algebra $R(S) = \underset{n=0}{\overset{}{}} R(S_n)$ is the central object rather then **Symm**.

## 2K. Product formulas

The simplest (arithmetic) product formula concerns the real and *p*-adic absolute values of a rational number

$$|a| = |a|_p^{-1}$$

where the product on the right is over all prime numbers $p$. There are more formulas of this type. This leads to a view of things that is expressed as follows by Yuri Manin in [129 Manin], Reflections on arithmetical physics, pp 149ff.

"Now we can see the following pattern

• (at least some) essential notions of real and complex calculus have their adèlic counterparts;

• adèlic objects have a strong tendency to be simpler than their Archimedean components, e.g. the adèlic fundamental domains of arithmetical discrete subgroups of semisimple groups usually have volume 1 (the Siegel-Tamagawa-Weil philosophy ... );

• due to this fact and to product formulas like (2) or (3) embodying the idea of democracy for all topologies, information on the real component of an adèlic object can be read off either from the real component or the product of the $p$-adic components for all $p$'s.

With some strain one can generalize and state the following principle which is the main conjecture of this talk.

*On the fundamental level our world is neither real, nor p-adic; it is adèlic. For some reasons reflecting the physical nature of our kind of living matter (e.g., the fact that we are built of massive particles), we tend to project the adèlic picture onto its real side. We can equally well spiritually project it upon its non-Archimedean side and calculate most important things arithmetically.*"

There are applications of this idea to the Polyakov measure (Polyakov partition function), loc. cit., string theory, [67 Freund et al.], Yang-Mills theory, [18 Asok et al.], and much more, see, for a start, (the bibliography of) [105 Khrennikov]. Add to this that the $p$-adic versions are often easier to handle and one finds some good justification for the discipline of $p$-adic physics.

## 3. Some first results and theorems.

### 3A. Freeness theorems .

The only general freeness theorem that I know about is the one from [66 Fressé]. This one says that cogroups (cogroup objects) in the category of algebras over an operad are free. This covers for instance one of the Kan results, the Leray theorem, the Milnor–Moore theorem and probably several more. At this stage it is unclear how far it goes.

I don't think it can be made to take care of the subobject freeness theorems; but there probably is a general theorem, yet to be formulated and proved, that can take care of those.

### 3B. On the Lazard universal formal group theorem.

The Lazard universal formal group theorem says that there exists a universal (one dimensional) commutative formal group (trivial) and that the underlying ring is free commutative polynomial in an infinity of indeterminates (surprising and far from trivial). The standard proof is long , laborious, and computational, even when simplified and streamlined as in [69 Fröhlich], see also [151 Ravenel].

Having a candidate universal formal group available, as in [80 Hazewinkel; 83 Hazewinkel] helps a great deal, see [80 Hazewinkel], pp 27 –30. But the proof is still mainly computational; also the construction of the candidate universal formal group involves choices of coefficients, which mars things. One dreams of a proof which mainly relies on universality properties.

In this connection there is a rather different proof due to Cristian Lenart, [118 Lenart], which seems to have promising aspects. One ingredient, which I consider promising, is the following. Consider the power series

$$f_b(X) = X + b_2 X^2 + b_3 X^3 + \cdots$$

over $\mathbf{Z}[b]$. Here the $b$'s are inderterminates. Now form

$$F_b(X,Y) = f_b(f_b^{-1}(X) + f_b^{-1}(Y))$$

This is of course a formal group over $\mathbf{Z}[b]$. It is proved [19] in loc. cit. that the coefficients of $F_b(X,Y)$ generate a free polynomial subring, $L$, of $\mathbf{Z}[b]$ and that regarded as a formal group over the subring $L$ $F_b(X,Y)$ is universal. Of course $L$ is truly smaller than $\mathbf{Z}[b]$. To start with $2b_2$ $L$, but $b_2$ $L$.

This next bit is pure speculation. The first Cartier theorem on formal groups says that the formal group of Witt vectors, $\hat{W}$, represents the functor 'curves'. This is a rather different universality property for formal groups. The covariant bialgebra of $\hat{W}$ is **Symm**. One wonders whether this can be used to prove the Lazard theorem.

### 3C. Objects and isomorphisms in connection with Symm.

This whole subsection is concerned with the objects and isomorphisms in the diagram at the end of section **2J**.

3C.1. *The isomorphism 'Ha' between* $R(W)$, *the representing ring of the functor of the big Witt vectors and* $U(\ )$, *the free lambda ring on one generator.*

Here is a synopsis of the relevant bits of structure. The ring $R(W)$ represents a (covariant) functor that carries a comonad structure, and the coalgebras for this comonad are precisely the lambda rings. That is all that is needed.

---

[19] The result is nice; I consider the proof highly unsatisfactory.

Let $\mathcal{C}$ be a category and let $(T, \mu, \varepsilon)$ be a comonad in $\mathcal{C}$. Now let $(Z, z \in T(Z))$ represent the functor $T$. That is, there is a functorial bijection
$\mathcal{C}(Z, A) \simeq T(A)$, $f \mapsto T(f)(z)$. The comonad structure gives in particular a morphism
$\sigma: Z \to TZ$, viz the image of $\mathrm{id}_Z$ under
$\mu_Z: T(Z) = \mathcal{C}(Z, Z) \to T(T(Z)) = \mathcal{C}(Z, T(Z))$. This defines a 'coalgebra for $T$' structure on $Z$. Now let $(A, \sigma)$ be a coalgebra for the comonad $T$ and let $a$ be an element of $A$. Consider the element $\sigma(a) \in T(A) = \mathcal{C}(Z, A)$. This gives a unique morphism of $T$-coalgebras that takes $z$ into $a$. There are of course a number of things to verify both at this categorical level and to check that these categorical considerations fit with the explicit constructions carried out in the previous subsections. This is straightforward.

Thus the isomorphism 'Ha' is a special case of a quite general theorem and the proof uses no special properties but only universal and other categorical notions. This is the kind of proof I would like to have for all the isomorphisms in the diagram.

3C.2. *The isomorphism 'Z' between $R(S)$ and* **Symm**. This is handled by the Zelevinsky theorem, [175 Zelevinsky] and [87 Hazewinkel], chapter 3. The Zelevinsky theorem deals with PSH algebras (over the integers). The acronym 'PSH' stands for 'Positive–Selfadjoint–Hopf. Actually it is about (nontrivial) graded connected positive self–adjoint Hopf algebras with a distinguished (preferred) homogenous basis. The Hopf algebra is also supposed to be of finite type so that each homogenous component is a free Abelian group of finite rank.

An inner product is defined by declaring this basis to be orthonormal. The positive elements of the Hopf algebra are the nonnegative (integer coefficient) linear combinations of the distinguished basis elements. Let $m$ and $\mu$ denote the multiplication and comultiplication respectively.

Selfadjoint (selfdual) now means

$$\langle m(x \otimes y), z \rangle = \langle x \otimes y, \mu(z) \rangle$$

and positivity means that if the elements of the distinguished basis are denoted by $\omega_i$ etc., and

$$m(\omega_i \otimes \omega_j) = \sum_r a_{i,j}^r \omega_r, \quad \mu(\omega_i) = \sum_{r,s} b_i^{r,s} \omega_r \otimes \omega_s$$

then $a_{i,j}^r \geq 0$ and $b_i^{r,s} \geq 0$.

Suppose now that there is precisely one among the distinguished basis elements that is primitive [20], then (the main part of) the Zelevinsky theorem says that the Hopf algebra in question is isomorphic (as a Hopf algebra) to **Symm**, possibly degree shifted.

_______________

[20] There is always at least one because of graded connectedness (and nontriviality).

An example of a PSH algebra is $R(S)$:

• The distinguished basis is formed by the irreducible representations of the various $S_n$

• The positive elements are the real (as opposed to virtual) representations, and so multiplication and comultiplication are positive.

• The selfadjointness comes from Frobenius reciprocity

• The Hopf property is handled by (a consequence of) the Mackey double coset theorem.

Using the isomorphism all structure can be transferred making **Symm** also a PSH algebra. An odd thing is that this is not proved directly. The distinguished basis turns out to be formed by the Schur functions. The problem is positivity. There seems to be no direct proof in the literature that the product of two Schur functions is a nonnegative linear combination of Schur functions.

I used to think that this theorem did not count in the context of the diagram because it uses such seemingly non-algebraic things as positivity and distinguished basis. However in the setting of $R(S)$ these are, see above, entirely natural.

There is one more thing I would like to say in this context. The fourth and final step of the proof of the Zelevinsky theorem (in the presentation of [87 Hazewinkel]) essential use is made of something called the Bernstein morphism. This is a morphism

$$H \qquad H \quad \textbf{Symm}$$

defined for any commutative associative graded connected Hopf algebra $H$. If one takes $H = \textbf{Symm}$ this is precisely the morphism that defines the multiplication on the big Witt vectors. This is a "coincidence" that cries out for further investigation.

For a completely different way of establishing that **Symm** and $R(S)$ are isomorphic see [19 Atiyah]. For still another and very elegant proof of this result see [121 Liulevicius; 122 Liulevicius]. It seems that the theorem actually goes back to Frobenius, [68 Frobenius].

3C.3. *The isomorphism 'S' from $R(S)$ to $R_{rat}(\textbf{GL})$*. This is Schur–Weyl duality which has its origins in Schur's thesis of 1901, [155 Schur]. The subject of Schur-Weyl duality has by now evolved into what is practically a small specialism of its own. A search in the MathSci database gave 72 hits. There are quantum and super versions and there are interrelations with such diverse fields as quantum and statistical mechanics, tilting theory, combinatorics, random walks on unitary groups, ... . A selection of references is [12; 26 Benkart et al.; 52 Dipper et al.; 54 Doty; 55 Duchesne; 70 Fulton et al.; 72 Goodman; 73 Goodman et al.; 74 Green; 93 Howe; 106 Klink et al.; 120 Lévy; 155 Schur; 154 Schur; 172 Weyl; 174 Zelditch], [19 Atiyah; 127 Macdonald; 59 Felder et. al.]

Here is what is probably the simplest incarnation of Schur–Weyl duality. Let $V$ be a finite dimensional vector space over a field of characteristic 0. Form the $n$-th tensor

product

$$T^n(V) = V \quad \cdots \quad V$$

The symmetric group $S_n$ acts on this by permuting the factors, which gives a finite dimensional representation of $S_n$ that can be decomposed into its isotypic components

$$T^n(V) = \quad \text{Hom}_{kS_n}(E_\pi, T^n(V)) \quad E_\pi \quad = \quad F_\pi(V) \quad E_\pi \qquad (3C3.1)$$
$$\pi \qquad\qquad\qquad\qquad\qquad\qquad\qquad\qquad \pi$$

functorially in $V$. Here the $E_\pi$ are the distinct irreducible $kS_n$ modules. If now $A : V \qquad V$ is a linear transformation $F_\pi(A)$: $F_\pi(V) \qquad F_\pi(V)$ is an 'invariant matrix' in the sense of Schur, [155 Schur]. This is taken from [127 Macdonald].

Taking invertible $A$ one obtains a representation $F_\pi(V)$ of $\mathbf{GL}(V)$. This can also be seen as coming from the action of $\mathbf{GL}(V)$ on $T^n(V)$ defined by $g(v_1 \quad \cdots \quad v_n) = gv_1 \quad \cdots \quad gv_n$, noting that this action commutes with the $S_n$ action on $T^n(V)$ and using the double commutant theorem.

The middle term in (3C.3.1) makes it clear that this is some kind of duality. What I would really like is to have is a pairing $R(S) \times R_{\text{rat}}(\mathbf{GL}) \qquad \mathbf{Z}$ defined directly, which then gives this duality. At the 'finite level' described above this can probably be done by looking at the trace form $\langle X, Y \rangle = \text{Trace}(XY)$ on $\text{End}(T^n(V))$, [73 Goodman et al.], section 9.1; [174 Zelditch], page 19. But not it seems without bringing in a lot of representation theory.

3C.4. *On a possible isomorphism 'L' between* $R(S)$ *and* $H(\mathbf{BU};\mathbf{Z})$. This is mostly speculative. First both rings (as Abelian groups) have a natural basis indexed by partitions. Second there is a bit of positive evidence in [79 Hazewinkel et al.], where in section 11 a (nontrivial, i.e. with jumps) family of representations is constructed of $S_{n+m}$ that is parametrized by the Grassmann manifold $\mathbf{Gr}_n(\mathbf{C}^{n+m})$.

3C.5. *On the isomorphism 'Du' between* $H(\mathbf{BU};\mathbf{Z})$ *and* $H(\mathbf{BU};\mathbf{Z})$. This is a matter of homology – cohomology duality for oriented manifolds. Plus autoduality of the Hopf algebras involved. (Both carry natural Hopf algebra structures.)

3C.6. *On the isomorphism 'SP' between* $H(\mathbf{BU};\mathbf{Z})$ *and* $\mathbf{Symm}$. First one defines Chern classes, for instance as in [138 Milnor et al.], chapter 14; see [8] for a totally different method; the definition of the first Chern class that is in [104 Kharshiladze] is one I particularly like.

The $i$-th Chern class associates to a complex vector bundle $V$ over a suitable space $X$ an element $c_i(V)$ of the cohomology group $H^{2i}(X;\mathbf{Z})$. One of the more important properties of the Chern classes is 'functoriality'. Let $f: Y \longrightarrow X$ be continuous and let $f^*V$ be the vector bundle pullback of $V$. Then

$$c_i(f^*V) = f^*(c_i(V))$$

(The notation is a bit unfortunate in that there are two different $f^*$ in the formula; but is traditional). A second important property is the 'Whitney sum formula'. Let

$$c(V) = 1 + c_1(V) + c_2(V) + \cdots$$

be the so-called total Chern class (also sometimes called complete Chern class). Let $W$ be a second complex vector bundle over $X$. Then

$$c(V \oplus W) = c(V)c(W)$$

where on the right hand side the cohomology cup product is used. And in fact together with $c_0(V) = 1$ and a normalization condition that specifies the total Chern class of the canonical (tautological) line bundle over the complex projective spaces $\mathbf{Gr}_1(\mathbf{C}^n)$ these two properties completely determine the Chern classes. See also [78 Hatcher], theorem 3.2 on page 78.

Next one calculates the cohomology of the classifying spaces $\mathbf{BU}_n$ to be

$$H^*(\mathbf{BU}_n;\mathbf{Z}) = \mathbf{Z}[c_1, c_2, \cdots, c_n]$$

where the $c_i$ are the Chern classes of the canonical vector bundle $\gamma_n$ over $\mathbf{BU}_n$. For instance with induction starting with the very simple case $\mathbf{BU}_1 = \mathbf{CP}^\infty$ which has a CW complex cell decomposition with precisely one cell in every even dimension. This is the way it is done in [138 Milnor et al.]. One can also use special sequences. It follows that

$$H^*(\mathbf{BU};\mathbf{Z}) = \mathbf{Z}[c_1, c_2, \cdots, c_n, \cdots]$$

which is isomorphic, at least as rings, to **Symm**. This is precisely the kind of calculatory proof that I do not like.

However, there is the following aspect. It is often a good idea to view **Symm** as the symmetric functions in an infinity of indeterminates

$$\mathbf{Symm} \subset \mathbf{Z}[\xi_1, \xi_2, \xi_3, \cdots]$$

Now on the topological side consider the canonical line bundle $\gamma_1$ $\mathbf{BU}_1$ and take the $n$-fold product $\gamma_1 \times \cdots \times \gamma_1$. This an $n$-dimensional bundle over the $n$-fold product $\mathbf{BU}_1 \times \cdots \times \mathbf{BU}_1$. The cohomology of this space is

$$\mathbf{Z}[\eta_1, \cdots, \eta_n]$$

where $\eta_i$ is the first Chern class of the $i$-th $\gamma_1$. Also by the Whitney sum formula

$$c(\gamma_1 \times \cdots \times \gamma_1) = (1 + \eta_1) \cdots (1 + \eta_n)$$

Now by the classifying space property of $\mathbf{BU}_n$ there is a homotopy class of maps $f : \mathbf{BU}_1 \times \cdots \times \mathbf{BU}_1$ $\mathbf{BU}_n$ such that the pullback of $\gamma_n$ by $f$ is $\gamma_1 \times \cdots \times \gamma_1$. Using functoriality it follows that $f$ takes $c_i$ $H$ $(\mathbf{BU}_n; \mathbf{Z})$ to the $i$-th elementary symmetric function in the $\eta$ 's and that $H$ $(\mathbf{BU}_n; \mathbf{Z})$ manifests itself as the ring of symmetric functions in $\mathbf{Z}[\eta_1, \cdots, \eta_n]$. This is taken from page 189 of [138 Milnor et al.]; see also [9] for a slightly different formulation of the same idea

Add to this that the Chern classes of the $\gamma_n$ (the universal Chern classes) can be described explicitly in terms of Schubert cycles, [8], and, possibly, this can be worked up to a much less calculatory proof of the isomorphism 'SP'.

3C.7. *The isomorphism 'F' between $R(W)$ and $U(\hat{W})$*. Consider an $n$-dimensional formal group $F$ over a (unital commutative associative) ring $A$. Here $n$ can be infinity. It is given by $n$ power series in $2n$ indeterminates grouped in two groups of $n$ indeterminates with coeficients in $A$. Let $R(F)$ be the ring of power series over $A$ in $n$ indeterminates. Then the $n$ power series of the formal group $F$ define a bialgebra like structure

$$R(F) \qquad R(F) \,\hat{}\, R(F)$$

This object is called the contravariant bialgebra of the formal group. (It is really needed (in general) to take the completed tensor product; even for $n = 1$ one has $A[[X]]$ $A[[Y]]$ $A[[X,Y]]$.)

$R(F)$ is given the usual power series topology. Now form

$$U(F) = \mathbf{Mod}_{A,\mathrm{cont}}(R(F), A) \qquad\qquad (3C.7.1)$$

This is the covariant bialgebra (in fact Hopf algebra) of the formal group $F$. Inversely one can obtain $R(F)$ from $U(F)$; just how will not be needed here.

In the case of the formal group $\hat{W}$ of the Witt vectors (over the integers) the power series defining it are in fact polynomials. And thus the restriction to

$$R(W) = \mathbf{Z}[X_1, X_2, \cdots] \qquad \mathbf{Z}[[X_1, X_2, \cdots]] = R(\hat{W})$$

of $R(\hat{W}) \qquad R(\hat{W}) \; \hat{} \; R(\hat{W})$ lands in $R(W) \quad R(W)$. As the polynomials are dense in the power series, in this polynomial case, formula (3C.7.1) is equivalent to

$$U(\hat{W}) = \mathbf{Mod}_{\mathbf{Z}}(R(W), \mathbf{Z})$$

and thus isomorphism 'F' is a consequence of the autoduality of $R(W) = \mathbf{Symm}$.

3C.8. *On the isomorphisms* 'M1' *and* 'M3' *between* $R(S)$, $K(\mathbf{P}_k)$ *and* $U(\ )$.
One sees from formula (3C.3.1) in subsection 3C.3 that each irreducible representation of $S_n$ defines a functor of $\mathbf{V}_k$ to itself that is polynomial. Here $\mathbf{V}_k$ is the category of finite dimensional vector spaces over the field $k$, and polynomial means that for each pair of vector spaces $U, V$ the mapping $F : \mathrm{Hom}(U, V) \qquad \mathrm{Hom}(F(U), F(V))$ is polynomial. Let now $\mathbf{P}_k$ be the category of polynomial functors $\mathbf{V}_k \qquad \mathbf{V}_k$ of bounded degree and $K(\mathbf{P}_k)$ its Grothendieck group. Then the remarks just made practically establish the isomorphism 'M1'.

Next, $K(\mathbf{P}_k)$ carries a $\lambda$-ring structure induced by composition with the exterior powers $\phantom{}^i : \mathbf{V}_k \qquad \mathbf{V}_k$. It turns out that it thus becomes the free $\lambda$-ring on one generator, [19 Atiyah; 128 Macdonald]. This is 'M3'.

It needs to be sorted out whether the composition of 'M1' and 'M3' equals the composition of 'Z' and 'Ha'.

The main aim of [128 Macdonald] is to generalize this in various ways. Let $A$ be a $k$ algebra, $\mathbf{V}_A$ the category of finitely generated projective left $A$ modules, $\mathbf{P}_A$ the category of polynomial functors $\mathbf{V}_A \qquad \mathbf{V}_k$ of bounded degree and $K(\mathbf{P}_A)$ its Grothendieck group. Then $K(\mathbf{P}_A)$ is the free $\lambda$-ring generated by the classes of the functors $P \mapsto E \quad_A P$ where $E$ runs through a complete set of non-isomorphic finite dimensional simple right $A$–modules.

When applied to the group ring of a finite group there is also the result that

$R(G \sim S_n)$ is the free $\lambda$-ring on the irreducible representation of $G$. (Here $G \sim S_n$ $\phantom{}_{n \; 0}$ is the wreath product of $G$ and $S_n$.

Thus 'M1' and 'M3' are just the implest cases of much more general results, which makes them nicer in my view.

3C.9. *On the object $E(\mathbf{Z})$ and the isomorphisms 'Ho1' and 'Ho2'*. Peter Hoffman noted that there is a nice functor $E$, denoted 'exp' in [91 Hoffman] that makes some of what went before more elegant.

Let **Ab** be the category of Abelian groups and **GrRing** that of (unital ungraded–commutative) graded rings. An object $R$ of **GrRing** is a direct sum of Abelian groups $R_i$ together with multiplications $R_i \times R_j \rightarrow R_{i+j}$ making $\oplus_i R_i$ a unital commutative ring. As in the case of the big Witt vectors one considers the "1-units". To be precise consider the functor

$$\wedge: \mathbf{GrRing} \rightarrow \mathbf{Ab} \text{ defined by } \hat{R} = 1 + \prod_{i=1} R_i \qquad (3C.9.1)$$

where the Abelian group structure is given by multiplication.

Note that the functor of the big Witt vectors is given by $S \mapsto S[[t]] \mapsto S[[t]]^{\wedge}$. What this means is completely unexplored.

The functor (3C.9.1) has a left adjoint $\mathbf{Ab} \rightarrow \mathbf{GrRing}$, here denoted $E$, so that there is the functorial equality

$$\mathbf{GrRing}(E(A), R) = \mathbf{Ab}(A, \hat{R})$$

As a left adjoint $E(A)$ should be thought of as some kind of free object and, as is so often the case with functors that are part of an adjunction it picks up all kinds of extra structure. In this case it is first of all a Hopf algebra (as happened with the universal enveloping algebra). This comes from the observation that

$$E(A \oplus B) = E(A) \otimes E(B) \text{ [21]}$$

$E(A)$ carries a natural $\lambda$-ring structure. (Though I find the construction very difficult and, frankly, definitely on the ugly side.) However it is worth exploring further as it goes through the notion of what the author calls an $\omega$-ring, a notion equivalent to that of a $\lambda$-ring but whose axioms only involve linear maps. This gives one a shot at solving a rather vexing matter. **Symm** is a $\lambda$-ring; it is also selfdual. So, morally speaking, there should be something like a 'dual $\lambda$-ring structure' on it.

Returning to the paper [91 Hoffman], the main theorems appear to be

$$\bigoplus_{n \geq 0} R(G \sim S_n) \cong E(R(G))$$

$$E(A) \text{ is the free } \lambda\text{-ring generated by } A$$

which are very nice results showing that the functor $E$ merits further attention.

---

[21] This formula also illustrates that 'exponential' or 'exp' is a most apt appellation.

3C.10. *The K-theory of endomorphisms.* Let $A$ be a unital commutative ring. Consider the category End(A) of pairs $(P, f)$ where $P$ is a finitely generated projective $A$-module and $f$ an endomorphism of $P$. A morphism $\varphi : (P, f) \to (Q, g)$ in End(A) is a morphism $\varphi$ of $A$-modules that commutes with the given endomorphisms, i.e. $g\varphi = \varphi f$. There is an obvious notion of exact sequence in End(A) and so one can form the Grothendieck group and ring [22] $K(\text{End}(A))$, the study of which was initiated by Gert Almkvist, [15 Almkvist; 16 Almkvist].

Given $(P, f) \in \text{End}(A)$ let $Q$ be a finitely generated module such that $P \oplus Q$ is free and consider the endomorphism $f \oplus 0$ of this module and its characteristic polynomial $\det(1 + t(f \oplus 0))$. This is a polynomial in $t$ that does not depend on $Q$. This induces a homomorphism $K(\text{End}(A)) \to W(A)$, where $W(-)$ is the functor of the big Witt vectors, that is obviously zero on $K(A)$. (The projective modules over $A$ are imbedded in End(A) as pairs $(A, 0)$). Thus there results a morphism (of rings in fact)

$$c : K(\text{End}(A))/K(A) = W_0(A) \to W(A)$$

functorial in $A$. Almkvist now proves:

The morphism $c$ is injective for all $A$ and the image of $c$ (for a given $A$) consists of all power series $1 + a_1 t + a_2 t^2 + \cdots$ that can be written in the form

$$1 + a_1 t + a_2 t^2 + \cdots = \frac{1 + b_1 t + b_2 t^2 + \cdots + b_r t^r}{1 + d_1 t + d_2 t^2 + \cdots + d_n t^n} \quad \text{with } b_i, d_j \in A$$

For obvious reasons I call these rational Witt vectors.

A first question is now whether this functor $W_0(-)$ is representable. It is, [85 Hazewinkel]. This requires some preparation. Consider the ring

$$\mathbf{Z}[X] = \mathbf{Z}[X_1, X_2, X_3, \cdots]$$

of polynomials in a countable infinity of commuting indeterminates. Form the Hankel matrix

$$\begin{matrix} 1 & X_1 & X_2 & X_3 & \cdots \\ X_1 & X_2 & X_3 & X_4 & \cdots \\ X_2 & X_3 & X_4 & X_5 & \cdots \\ \vdots & \vdots & \vdots & \vdots & \ddots \end{matrix}$$

---

[22] The multiplication is induced by the tensor product.

Now let $J_n$ be the ideal in $\mathbf{Z}[X]$ generated by all $(n+1) \times (n+1)$ minors of this Hankel matrix. These ideals define a topology on $\mathbf{Z}[X]$ which for the present purposes I will call the $J$-topology. The representability result is now as follows.

For each rational Witt vector $a(t) = 1 + a_1t + a_2t^2 + \cdots$ $W_0(A)$ let

$\varphi_{a(t)} : \mathbf{Z}[X]$ $A$ be the ring morphism defined by $X_i \mapsto a_i$. Then $a(t) \mapsto \varphi_{a(t)}$ is a functorial and injective morphism from $W_0(A)$ to ring morphisms $\mathbf{Z}[X]$ $A$ that are continuous with respect to the $J$-topology on $\mathbf{Z}[X]$ and the discrete topology on $A$. If $A$ is Fatou, so in particular if $A$ is integral and Noetherian, the correspondence is bijective.

Here Fatou is a technical condition that is of no particular importance for this paper. Suffice it to say that a Noetherian integral domain is Fatou. Incidentally the quotient rings $\mathbf{Z}[X]/J_n$ are integral domains, but they are not Noetherian and not Fatou.

For a host of other results, including a determination of the operations in the $K$-theory of endomorphisms, see [15 Almkvist; 16 Almkvist; 85 Hazewinkel].

3C.11. *Leftovers.*

• **Symm** is an object with an enormous amount of compatible structure: Hopf algebra, inner product, selfdual (as a Hopf algebra), PSH, coring object in the category of rings, ring object in the category of corings (up to a little bit of unit trouble), Frobenius and Verschiebung endomorphisms, free algebra on the cofree coalgebra over $\mathbf{Z}$ (and the dual of this: cofree coalgebra over the free algebra on one element), several levels of lambda ring structure, ... .

The question arises which ones of these have natural interpretations in the other nine incarnations occurring in the diagram (and whether the isomorphisms indicated are the right ones for preserving these structures).

• **Symm** represents the functor of the big Witt vectors

$W(A) = \{1 + a_1t + a_2t^2 + \cdots : a_i \quad A\}$.

Now **Hopf**(**Symm**, **Symm**) $= W(\mathbf{Z})$, [121 Liulevicius]. This comes about because on the one hand **Symm** is the free algebra on the cofree coalgebra over $\mathbf{Z}$, and on the other the cofree coalgebra over the free algebra over $\mathbf{Z}$.

This is a curiosity that certainly merits some thought and one wonders whether something similar occurs elsewhere.

## References

The list of references below contains more items than are actually referred to in the text above. The others are included because I know or suspect that there are more niceness results in them.